\documentclass[twoside,english,3p]{elsarticle}
\usepackage[T1]{fontenc}
\usepackage{geometry}
\usepackage{amssymb}
\usepackage{amsmath}

\usepackage{amsthm}
\usepackage{stmaryrd}
\usepackage{graphicx}
\usepackage{booktabs}
\usepackage{multirow}
\usepackage{makecell}
\usepackage{xurl}
\usepackage{esint}
\usepackage{algorithm}
\usepackage{algpseudocode}
\usepackage{rotating}
\usepackage{adjustbox}
\usepackage{nomencl}
\usepackage{multicol}
\usepackage{listings}
\makenomenclature
\makeatletter
\theoremstyle{plain}

\theoremstyle{boldremark} 

\ifx\proof\undefined

\providecommand{\proofname}{Proof}
\fi

\renewenvironment{proof}[1][\proofname]{%
	\par\pushQED{\qed}\normalfont%
	\topsep6\p@\@plus6\p@\relax
	\trivlist\item[\hskip\labelsep\bfseries#1\@addpunct{.}]%
	\ignorespaces
}{%
	\popQED\endtrivlist\@endpefalse
}

\journal{Elsevier}

\usepackage{hyperref}
\hypersetup{colorlinks = true, allcolors = blue}

\usepackage[nameinlink]{cleveref}

\Crefname{figure}{Fig.}{Figs.}
\Crefformat{equation}{Eq.~#2(#1)#3}
\Crefformat{section}{Section~#2#1#3}
\AtBeginDocument{%
	\let\citet\cite
}

\usepackage{float}
\newfloat{listing}{htbp}{lol}
\floatname{listing}{Listing}
\crefname{listing}{Listing}{Listings}
\Crefname{listing}{Listing}{Listings}

\usepackage[labelfont=bf]{caption}
\@ifundefined{showcaptionsetup}{}{%
	\PassOptionsToPackage{caption=false}{subfig}}
\usepackage{subfig}
\makeatother

\usepackage{babel}
\providecommand{\remarkname}{Remark}
\providecommand{\theoremname}{Theorem}

\begin{document}
	
	\begin{frontmatter}{}
		
		\title{Deep Energy Method with Large Language Model assistance: an open-source Streamlit-based platform for solving variational PDEs}

		\author[rvt,rvt3]{Yizheng Wang}
		
		\ead{wang-yz19@tsinghua.org.cn}
		
		\author[rvt3]{Cosmin Anitescu}
		
		\author[rvt6]{Mohammad Sadegh Eshaghi}
		
		\author[rvt6]{Xiaoying Zhuang}
		
		\author[rvt3]{Timon Rabczuk}
		\ead{timon.rabczuk@uni-weimar.de}
		
		\author[rvt]{Yinghua Liu\corref{cor1}}
		
		\ead{yhliu@mail.tsinghua.edu.cn}
		\cortext[cor1]{Corresponding author}
		\address[rvt]{Department of Engineering Mechanics, Tsinghua University, Beijing 100084, China}

		\address[rvt3]{Institute of Structural Mechanics, Bauhaus-Universit\"{a}t Weimar, Marienstr. 15, D-99423 Weimar, Germany}

		\address[rvt6]{ Institute of Photonics, Department of Mathematics and Physics, Leibniz University Hannover, Germany}

		\begin{abstract}
			Physics-informed neural networks (PINNs) in energy form, also known as the deep energy method (DEM), offer advantages over strong-form PINNs such as lower-order derivatives and fewer hyperparameters, yet dedicated and user-friendly software for energy-form PINNs remains scarce. To address this gap, we present \textbf{LM-DEM} (Large-Model-assisted Deep Energy Method), an open-source, Streamlit-based platform for solving variational partial differential equations (PDEs) in computational mechanics. LM-DEM integrates large language models (LLMs) for geometry modeling: users can generate Gmsh-compatible geometries directly from natural language descriptions or images, significantly reducing the burden of traditional geometry preprocessing. The solution process is driven by the deep energy method, while finite element solutions can be obtained in parallel. The framework supports built-in problems including Poisson, screened Poisson, linear elasticity, and hyperelasticity in two and three dimensions, as well as user-defined energy functionals analogous to the \texttt{UMAT} interface in Abaqus. The source code is available at \url{https://github.com/yizheng-wang/LMDEM}, and a web-based version is accessible at \url{https://ai4m.llmdem.com}. LM-DEM aims to lower the barrier for practitioners and beginners to adopt energy-form PINNs for variational PDE problems.
		\end{abstract}

		\printnomenclature

		
	\end{frontmatter}{}
	
	\section{Introduction}
	
	A wide range of physical phenomena are governed by partial differential equations (PDEs), and the numerical solution of PDEs lies at the core of computational mechanics, playing a crucial role in understanding and predicting physical behaviors. Among existing numerical approaches, the finite element method (FEM) is the most widely used technique due to its high accuracy, strong numerical robustness, and solid theoretical foundation.
	
	Recently, AI for PDEs, one of the most important directions in AI for science, has emerged as a new paradigm that leverages deep learning techniques to solve PDEs. Broadly speaking, AI for PDEs can be categorized into three major classes: physics-informed neural networks (PINNs) \citep{PINN_original_paper,loss_is_minimum_potential_energy}, operator learning methods \citep{DeepOnet,li2020fourier}, and physics-informed neural operators (PINOs) \citep{li2024physics,eshaghi2025variational}. 
	Recently, hybrid approaches have emerged that first train a physics-informed neural operator (PINO) and then use its prediction as an initial guess for subsequent refinement with an iterative finite element solver, such as NOWS \citep{eshaghi2025nows} and P-FEM \citep{wang2026pretrain}.
	
	Among these approaches, PINNs have received the most attention and have become one of the most influential methodologies in AI for PDEs. PINNs approximate the solution of PDEs using neural networks and enforce physical laws through loss functions derived from governing equations and boundary conditions.
	PINNs can be further divided into two main formulations: the strong form \citep{PINN_original_paper} and the energy (variational) form, commonly referred to as the deep energy method (DEM) \citep{loss_is_minimum_potential_energy}. Given the significant interest in PINNs, the development of PINN-related software and libraries is of great importance, especially to facilitate rapid adoption by beginners and practitioners. In terms of existing PINN libraries and software, Lu et al.~\citep{PINNlibrary} developed DeepXDE based on TensorFlow, Haghighat et al.~\citep{haghighat2021sciann} proposed SciANN built upon Keras, and Zubov et al.~\citep{zubov2021neuralpde} constructed NeuralPDE in the Julia ecosystem. In addition, Hennigh et al.~\citep{hennigh2021nvidia} developed SimNet, which targets not only academic research but also industrial applications. 
	JAX provides a NumPy-like interface with automatic differentiation, making it highly suitable for differentiable PDE solvers.
	These libraries significantly lower the barrier for applying PINNs to solve PDEs.
	
	Moreover, numerous open-source implementations of PINNs are available, as summarized in Table~\ref{tab:AI4PDEcodes}. However, it can be observed that the vast majority of existing software and libraries focus on the strong-form PINNs. In contrast, there is a notable lack of dedicated libraries and user-friendly software for energy-form PINNs, especially those designed for ease of use by beginners. Compared with strong-form PINNs, energy-form PINNs typically involve lower-order derivatives and require substantially fewer hyperparameters \citep{the_comparision_of_strong_and_energy_form,wang2022cenn}, which often leads to improved numerical accuracy, stability, and computational efficiency. Therefore, the development of software specifically tailored for energy-form PINNs is both necessary and timely.
	
	\begin{table}
		\centering{}\caption{The codes of PINNs strong and energy form in computational mechanics\label{tab:AI4PDEcodes}}
		\begin{adjustbox}{max width=\textwidth}
			\begin{tabular}{ccc}
				\toprule 
				& The link of code & Brief Description of Method\tabularnewline
				\midrule
				\multirow{11}{*}{Strong form} & \url{https://github .com /maziarraissi /PINNs} & The original code\citep{PINN_original_paper}\tabularnewline
				& \url{https://github.com/ehsankharazmi/hp-VPINNs} & Variational PINNs\citep{hp-VPINN}\tabularnewline
				& \url{https://github.com/AmeyaJagtap/Conservative_PINNs} & PINNs with subdomains\citep{CPINN}\tabularnewline
				& \url{https://github.com/lululxvi/hpinn} & Inverse design\citep{lu2021physics}\tabularnewline
				& \url{https://github.com/XuhuiM/PPINN} & Time-dependent PDEs\citep{meng2020ppinn}\tabularnewline
				& \url{https://github.com/sciann/sciann-applications} & Solid mechanics\citep{PINN_solid_mechanics}\tabularnewline
				& \url{https://github.com/Raocp/PINN-elastodynamics } & Elastodynamics\citep{PINNstrong_form_in_elastodynamics}\tabularnewline
				& \url{https://github.com/Jianxun-Wang/LabelFree-DNN-Surrogate } & Incompressible NS \citep{sun2020surrogate}\tabularnewline
				& \url{https://github.com/PredictiveIntelligenceLab/DeepStefan } & Multiphase and moving boundary problems\citep{wang2020deep}\tabularnewline
				& \url{https://github.com/Raocp/PINN-laminar-flow} & Incompressible laminar flows\citep{rao2020physics}\tabularnewline
				& \url{https://github.com/shengzesnail/AIV_MAOAC} & Blood flow\citep{cai2021artificial}\tabularnewline
				\cmidrule{2-3} \cmidrule{3-3} 
				\multirow{9}{*}{Energy form} & \url{https://github.com/ISM-Weimar/DeepEnergyMethods} & The original code\citep{loss_is_minimum_potential_energy}\tabularnewline
				& \url{https://github.com/yizheng-wang} & DCEM\citep{wang2023dcm}\tabularnewline
				& \url{https://github.com/yizheng-wang} & DEM with subdomains\citep{wang2022cenn}\tabularnewline
				& \url{https://github.com/weili101/Deep Plates} & Kirchhoff plate\citep{zhuang2021deep,the_comparision_of_strong_and_energy_form}\tabularnewline
				& \url{https://github.com/somdattagoswami/IGAPack-PhaseField} & Phase field for fracture mechanics\citep{goswami2020adaptive,goswami2020transfer}\tabularnewline
				& \url{https://github.com/MinhNguyenIKM/dem_hyperelasticity} & Hyperelasticity\citep{PINN_hyperelasticity}\tabularnewline
				& \url{https://github.com/JinshuaiBai/RPIM_NNS} & Hyperelasticity\citep{bai2024robust}\tabularnewline
				& \url{https://github.com/Jasiuk-Research-Group/DEM_for_J2_plasticity} & J2 elastoplasticity\citep{he2023deep}\tabularnewline
				& \url{https://github.com/MinhNguyenIKM/parametric-deep-energy-method} & Parametric DEM\citep{paradeepenergy}\tabularnewline
				&
				\url{https://github.com/eshaghi-ms/DeepNetBeam} & Functionally graded porous beams \citep{eshaghi2025applications}  \tabularnewline
				\bottomrule
			\end{tabular}
		\end{adjustbox}
	\end{table}
	
	To address the aforementioned gap, we develop an open-source computational mechanics software based on energy-form PINNs, termed \textbf{LM-DEM} (Large-Model-assisted Deep Energy Method). LM-DEM is an AI-driven framework in which geometry modeling is performed using large language models, while the solution process is based on the deep energy method. In addition, LM-DEM is capable of producing finite element solutions for comparison and validation. The open-source code of LM-DEM is available at \url{https://github.com/yizheng-wang/LMDEM}, and the web-based version can be accessed at \url{https://ai4m.llmdem.com}.
	
	The main features of LM-DEM are summarized as follows:
	\begin{itemize}
		\item \textbf{Large-model-assisted geometry modeling}: both 2D and 3D geometries can be generated directly from natural language descriptions or images.
		\item \textbf{AI-driven solution}: PDEs are solved using the deep energy method, while finite element solutions can also be obtained for reference.
		\item \textbf{Customizable energy functionals}: users can define and modify energy functionals according to their specific physical problems.
	\end{itemize}
	
	The remainder of this paper is organized as follows. Section~\ref{sec:Preparatory-knowledge} introduces the strong-form and energy-form PINNs. Section~\ref{sec:Method} presents a user guide for LM-DEM. Section~\ref{sec:Result} is divided into three parts: (i) geometry generation examples using large language models, (ii) numerical solutions obtained by the deep energy method and the finite element method, and (iii) the usage of custom energy functionals together with corresponding numerical results. Finally, conclusions are drawn in Section~\ref{sec:Conclusion}.

	\section{Preparatory Knowledge}
	\label{sec:Preparatory-knowledge}
	
	In this section, we briefly review the strong-form and energy-form formulations of physics-informed neural networks (PINNs), with particular emphasis on the energy form. This focus is motivated by the fact that LM-DEM constructs its loss function based on variational principles and energy minimization.
	
	\subsection{Strong Form}
	
	We begin with a general boundary value problem governed by partial differential equations:
	\begin{equation}
		\begin{cases}
			\boldsymbol{L}(\boldsymbol{u}(\boldsymbol{x}))=\boldsymbol{f}(\boldsymbol{x}), & \boldsymbol{x}\in\Omega,\\
			\boldsymbol{B}(\boldsymbol{u}(\boldsymbol{x}))=\boldsymbol{g}(\boldsymbol{x}), & \boldsymbol{x}\in\Gamma,
		\end{cases}
		\label{eq:original_form}
	\end{equation}
	where $\boldsymbol{L}$ and $\boldsymbol{B}$ denote the differential operators in the domain and on the boundary, respectively, while $\Omega$ and $\Gamma$ represent the computational domain and its boundary.
	
	Using the weighted residual method, Eq.~\eqref{eq:original_form} can be rewritten as
	\begin{equation}
		\begin{cases}
			\int_{\Omega}\big[\boldsymbol{L}(\boldsymbol{u}(\boldsymbol{x}))-\boldsymbol{f}(\boldsymbol{x})\big]\cdot\boldsymbol{w}(\boldsymbol{x})\,\mathrm{d}\Omega=0, & \boldsymbol{x}\in\Omega,\\
			\int_{\Gamma}\big[\boldsymbol{B}(\boldsymbol{u}(\boldsymbol{x}))-\boldsymbol{g}(\boldsymbol{x})\big]\cdot\boldsymbol{w}(\boldsymbol{x})\,\mathrm{d}\Gamma=0, & \boldsymbol{x}\in\Gamma,
		\end{cases}
		\label{eq:original_form_weighted}
	\end{equation}
	where $\boldsymbol{w}(\boldsymbol{x})$ is an arbitrary weighting function. When $\boldsymbol{w}$ is chosen arbitrarily, the weighted residual form is mathematically equivalent to the original PDE system.
	
	For numerical convenience, a specific choice of the weighting function is commonly adopted:
	\begin{equation}
		\boldsymbol{w}(\boldsymbol{x})=
		\begin{cases}
			\boldsymbol{L}(\boldsymbol{u}(\boldsymbol{x}))-\boldsymbol{f}(\boldsymbol{x}), & \boldsymbol{x}\in\Omega,\\
			\boldsymbol{B}(\boldsymbol{u}(\boldsymbol{x}))-\boldsymbol{g}(\boldsymbol{x}), & \boldsymbol{x}\in\Gamma.
		\end{cases}
		\label{eq:strong_form_weight}
	\end{equation}
	Substituting Eq.~\eqref{eq:strong_form_weight} into Eq.~\eqref{eq:original_form_weighted} yields
	\begin{equation}
		\begin{cases}
			\int_{\Omega}\big|\boldsymbol{L}(\boldsymbol{u}(\boldsymbol{x}))-\boldsymbol{f}(\boldsymbol{x})\big|^{2}\,\mathrm{d}\Omega=0,\\
			\int_{\Gamma}\big|\boldsymbol{B}(\boldsymbol{u}(\boldsymbol{x}))-\boldsymbol{g}(\boldsymbol{x})\big|^{2}\,\mathrm{d}\Gamma=0.
		\end{cases}
		\label{eq:original_form_delta}
	\end{equation}
	
	Approximating the integrals in Eq.~\eqref{eq:original_form_delta} via numerical quadrature or sampling leads to the loss function of the strong-form PINNs:
	\begin{equation}
		\mathcal{L}_{\mathrm{PINNs}}
		=
		\frac{\lambda_{r}}{N_{r}}\sum_{i=1}^{N_{r}}
		\big|\boldsymbol{L}(\boldsymbol{u}(\boldsymbol{x}_{i};\boldsymbol{\theta}))-\boldsymbol{f}(\boldsymbol{x}_{i})\big|^{2}
		+
		\frac{\lambda_{b}}{N_{b}}\sum_{i=1}^{N_{b}}
		\big|\boldsymbol{B}(\boldsymbol{u}(\boldsymbol{x}_{i};\boldsymbol{\theta}))-\boldsymbol{g}(\boldsymbol{x}_{i})\big|^{2}.
	\end{equation}
	The neural network approximation of the field variable $\boldsymbol{u}(\boldsymbol{x})$ is then obtained by
	\begin{equation}
		\boldsymbol{u}(\boldsymbol{x};\boldsymbol{\theta})
		=
		\underset{\boldsymbol{\theta}}{\arg\min}\;
		\mathcal{L}_{\mathrm{PINNs}}.
	\end{equation}
	
	From a mathematical perspective, the strong form of PINNs \citep{PINN_original_paper} corresponds to a specific choice of the weighting function in the weighted residual formulation, namely the residual of the governing equations themselves, as defined in Eq.~\eqref{eq:strong_form_weight}.
	
	\subsection{Energy Form}
	
	The strong-form PINNs suffer from two major drawbacks: the presence of multiple penalty hyperparameters and the requirement of higher-order derivatives of the neural network outputs. These issues can be alleviated by reformulating the problem using variational principles, leading to the energy form, which constitutes the foundation of the deep energy method (DEM).
	
	By choosing the weighting function in Eq.~\eqref{eq:original_form_weighted} as $\boldsymbol{w}=\delta\boldsymbol{u}$, the Galerkin formulation is obtained:
	\begin{equation}
		\int_{\Omega}\big[\boldsymbol{L}(\boldsymbol{u}(\boldsymbol{x}))-\boldsymbol{f}(\boldsymbol{x})\big]\cdot\delta\boldsymbol{u}\,\mathrm{d}\Omega=0,
		\quad \boldsymbol{x}\in\Omega.
		\label{eq:galerkin_form}
	\end{equation}
	
	To illustrate the idea, we consider the Poisson equation:
	\begin{equation}
		\begin{cases}
			-\Delta u(\boldsymbol{x})=f(\boldsymbol{x}), & \boldsymbol{x}\in\Omega,\\
			u(\boldsymbol{x})=\bar{u}(\boldsymbol{x}), & \boldsymbol{x}\in\Gamma^{u},\\
			\dfrac{\partial u(\boldsymbol{x})}{\partial n}=\bar{t}(\boldsymbol{x}), & \boldsymbol{x}\in\Gamma^{t},
		\end{cases}
		\label{eq:poisson_equation}
	\end{equation}
	where $\Gamma^{u}$ and $\Gamma^{t}$ denote the Dirichlet and Neumann boundaries, respectively.
	
	The Galerkin form of Eq.~\eqref{eq:poisson_equation} reads
	\begin{equation}
		\int_{\Omega}\big[-\Delta u(\boldsymbol{x})-f(\boldsymbol{x})\big]\delta u\,\mathrm{d}\Omega=0.
		\label{eq:galerkin_form_poisson}
	\end{equation}
	Applying integration by parts yields
	\begin{equation}
		\int_{\Omega}u_{,i}\,\delta u_{,i}\,\mathrm{d}\Omega
		-
		\int_{\Gamma}u_{,i}n_{i}\,\delta u\,\mathrm{d}\Gamma
		-
		\int_{\Omega}f\,\delta u\,\mathrm{d}\Omega
		=0.
		\label{eq:gaussian_poisson}
	\end{equation}
	By enforcing the Neumann boundary condition on $\Gamma^{t}$ and noting that $\delta u=0$ on $\Gamma^{u}$, we obtain the weak form
	\begin{equation}
		\int_{\Omega}u_{,i}\,\delta u_{,i}\,\mathrm{d}\Omega
		-
		\int_{\Gamma^{t}}\bar{t}\,\delta u\,\mathrm{d}\Gamma
		-
		\int_{\Omega}f\,\delta u\,\mathrm{d}\Omega
		=0.
		\label{eq:weak_form}
	\end{equation}
	
	Since $u(\boldsymbol{x})$ is prescribed on $\Gamma^{u}$, the admissible trial function must satisfy the Dirichlet boundary condition a priori. This requirement is essential in DEM. Moreover, Eq.~\eqref{eq:weak_form} naturally incorporates the domain equations and Neumann boundary conditions, making it equivalent to the original boundary value problem.
	
	The weak form can be equivalently derived from the variational principle
	\begin{equation}
		\delta\mathcal{L}
		=
		\int_{\Omega}u_{,i}\,\delta u_{,i}\,\mathrm{d}\Omega
		-
		\int_{\Gamma^{t}}\bar{t}\,\delta u\,\mathrm{d}\Gamma
		-
		\int_{\Omega}f\,\delta u\,\mathrm{d}\Omega,
		\label{eq:first_vari}
	\end{equation}
	where the potential energy functional is given by
	\begin{equation}
		\mathcal{L}
		=
		\frac{1}{2}\int_{\Omega}u_{,i}u_{,i}\,\mathrm{d}\Omega
		-
		\int_{\Gamma^{t}}\bar{t}u\,\mathrm{d}\Gamma
		-
		\int_{\Omega}fu\,\mathrm{d}\Omega.
		\label{eq:energy}
	\end{equation}
	We note that this equivalence between the Galerkin weak form and a variational (stationarity) principle holds specifically for problems that derive from an energy functional. For non-variational PDEs, a weak form can still be formulated (e.g., via weighted residual methods), but it does not, in general, correspond to the stationarity of an energy functional.
	
	For the Poisson equation considered here, since $\delta^{2}\mathcal{L}>0$, the solution can be obtained by minimizing the energy functional:
	\begin{equation}
		u(\boldsymbol{x})
		=
		\underset{u}{\arg\min}\;\mathcal{L}.
		\label{eq:minimum_energy}
	\end{equation}
	
	In DEM, the solution is approximated by a neural network $u(\boldsymbol{x};\boldsymbol{\theta})$, and the optimization problem becomes
	\begin{equation}
		\begin{aligned}
			u(\boldsymbol{x};\boldsymbol{\theta})
			&=
			\underset{\boldsymbol{\theta}}{\arg\min}\;\mathcal{L}_{\mathrm{DEM}},\\
			\mathcal{L}_{\mathrm{DEM}}
			&=
			\frac{1}{2}\int_{\Omega}u_{,i}(\boldsymbol{x};\boldsymbol{\theta})u_{,i}(\boldsymbol{x};\boldsymbol{\theta})\,\mathrm{d}\Omega
			-
			\int_{\Gamma^{t}}\bar{t}u(\boldsymbol{x};\boldsymbol{\theta})\,\mathrm{d}\Gamma
			-
			\int_{\Omega}f\,u(\boldsymbol{x};\boldsymbol{\theta})\,\mathrm{d}\Omega.
		\end{aligned}
		\label{eq:DEM}
	\end{equation}
	
	The core of DEM thus lies in the accurate evaluation of domain and boundary energies, as well as the construction of admissible trial functions. Numerical integration can be performed using Monte Carlo sampling or higher-order quadrature rules such as Gaussian integration or Simpson’s rule.
	
	To enforce the Dirichlet boundary condition, we adopt an admissible function based on the distance function:
	\begin{equation}
		u(\boldsymbol{x})
		=
		u_{p}(\boldsymbol{x};\boldsymbol{\theta}_{p})
		+
		D(\boldsymbol{x})\,u_{g}(\boldsymbol{x};\boldsymbol{\theta}_{g}),
		\label{eq:admissible}
	\end{equation}
	where $u_{p}$ is a particular solution network satisfying the Dirichlet boundary condition, trained via
	\begin{equation}
		\boldsymbol{\theta}_{p}
		=
		\underset{\boldsymbol{\theta}_{p}}{\arg\min}\;
		\mathrm{MSE}\big(u_{p}(\boldsymbol{x};\boldsymbol{\theta}_{p}),\bar{u}(\boldsymbol{x})\big),
		\quad \boldsymbol{x}\in\Gamma^{u}.
	\end{equation}
	The distance function $D(\boldsymbol{x})$ is approximated using radial basis functions \citep{wang2022cenn}:
	\begin{equation}
		D(\boldsymbol{x})
		=
		\min_{\boldsymbol{y}\in\Gamma^{u}}
		\|\boldsymbol{x}-\boldsymbol{y}\|.
		\label{eq:distance_function}
	\end{equation}
	The network $u_{g}(\boldsymbol{x};\boldsymbol{\theta}_{g})$ is a standard neural network, and only $\boldsymbol{\theta}_{g}$ is optimized in the energy minimization process:
	\begin{equation}
		u(\boldsymbol{x};\boldsymbol{\theta}_{p},\boldsymbol{\theta}_{g})
		=
		\underset{\boldsymbol{\theta}_{g}}{\arg\min}\;
		\mathcal{L}_{\mathrm{DEM}}.
	\end{equation}

	\section{Deep Energy Method with Large Language Model}
	\label{sec:Method}
	
\begin{figure}
		\centering
		\includegraphics[scale=0.6]{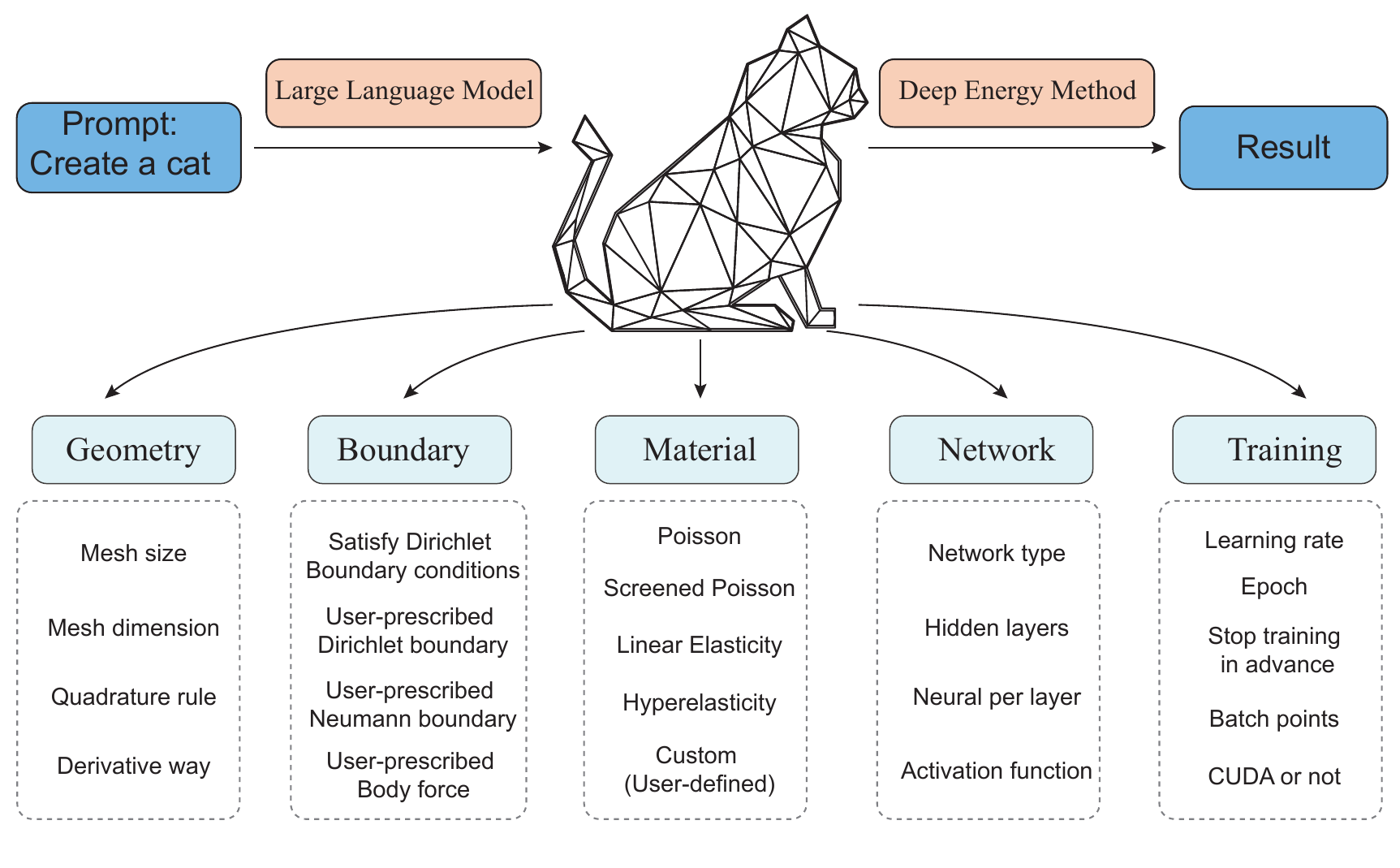}
		\caption{Overall framework of LM-DEM. The user provides a natural-language prompt (e.g., ``Create a cat'') to the large language model, which returns a geometry. Configuration is organized into five panels: Geometry (mesh dimension, quadrature rule, derivative evaluation, mesh size), Material (Poisson, screened Poisson, linear elasticity, hyperelasticity, or user-defined), Boundary (Dirichlet, Neumann, and body force), Network (layers, width, activation, type), and Training (epochs, early stopping, batch size, CUDA, learning rate). Dirichlet boundary conditions are enforced a priori; the pipeline then feeds into the deep energy method to obtain the solution.}
		\label{fig:LM-DEM}
	\end{figure}

	This section presents the practical usage of the proposed Deep Energy Method with Large Language Model assistance (LM-DEM). The workflow can be broadly divided into three stages: pre-processing, solver configuration, and post-processing. Figure~\ref{fig:LM-DEM} illustrates the overall software architecture of LM-DEM.
	
	LM-DEM can be used in two ways. First, users can deploy the open-source code from \url{https://github.com/yizheng-wang/LMDEM} on a local machine, which enables GPU acceleration. Second, users can directly access the web version at \url{https://ai4m.llmdem.com}; however, the web service currently supports CPU computation only, resulting in reduced efficiency. The local deployment is straightforward and can be launched by
	\begin{lstlisting}
		streamlit run DEM.py
	\end{lstlisting}
	
	\subsection{Pre-processing: Large Language Model-assisted Geometry Modeling}
	
	Geometry modeling is widely recognized as one of the most time-consuming components in computational simulation pipelines, often requiring extensive user interactions through numerous buttons and manual operations. To alleviate this burden, LM-DEM integrates large language model (LLM) assistance for geometry creation. Specifically, the user first selects the problem dimension (e.g., 2D) in the \texttt{Geometry} panel and then provides either a natural-language prompt or an input image. LM-DEM subsequently generates the target geometry automatically.
	
\begin{figure}
		\centering
		\includegraphics[scale=0.55]{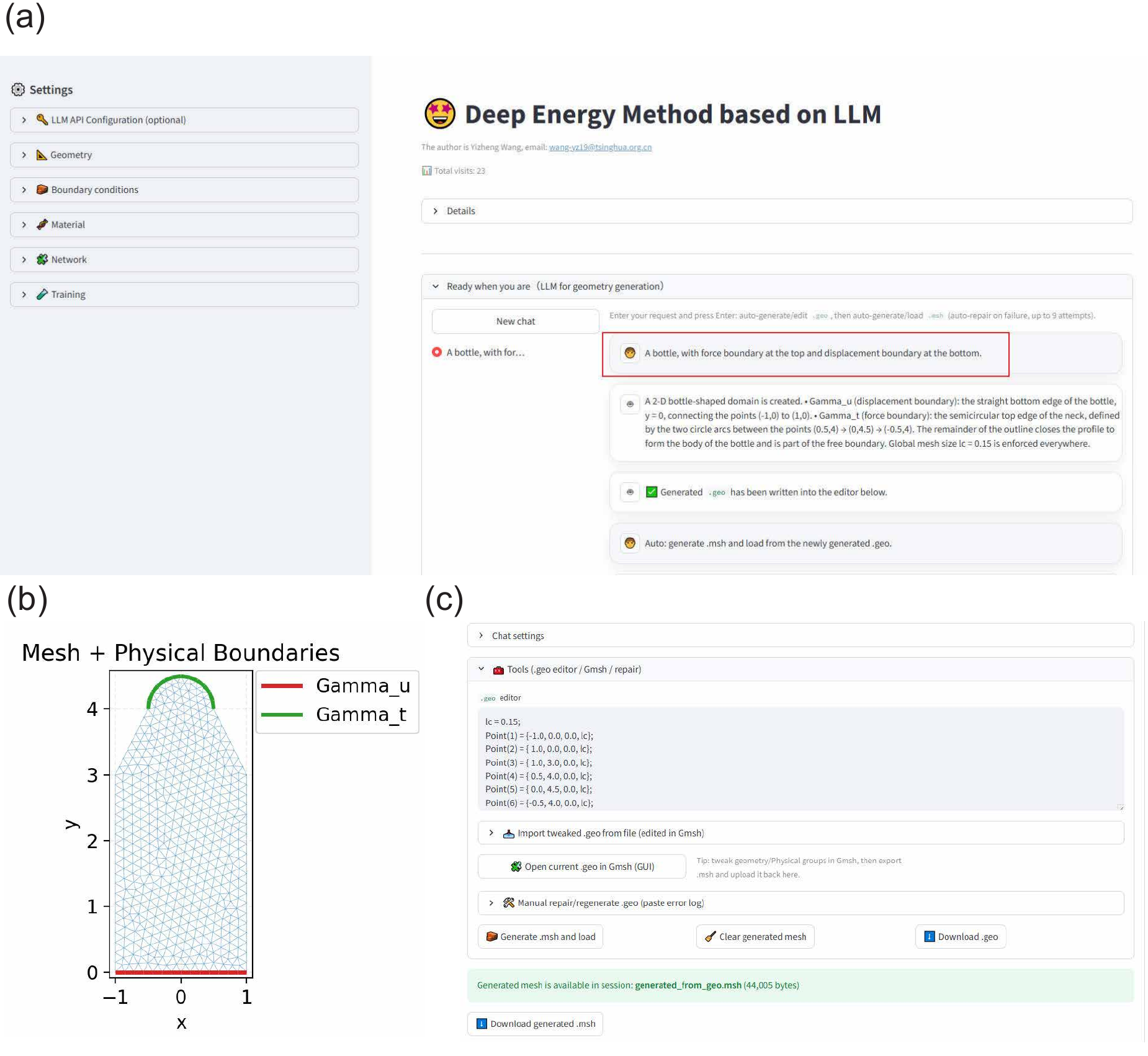}
		\caption{LLM-assisted geometry modeling in LM-DEM: (a) the main interface where users provide prompts specifying geometry and boundary-condition locations; the LLM returns a Gmsh-compatible \texttt{.geo} script, and users can switch backends (e.g., OpenAI, DeepSeek, Ollama) or start a new chat to retry; (b) the geometry generated and meshed via Gmsh; (c) the generated \texttt{.geo} script, which can be edited directly or opened in the Gmsh GUI for refinement.}
		\label{fig:LLM_generate}
	\end{figure}

	As an example, the prompt
	\begin{quote}
		\small ``A bottle, with force boundary at the top and displacement boundary at the bottom.''
	\end{quote}
	leads to the geometry generation process illustrated in Figure~\ref{fig:LLM_generate}. The core mechanism is to prompt an LLM to produce a \texttt{.geo} script that can be opened by Gmsh (\url{https://gmsh.info}). LM-DEM calls the corresponding LLM API automatically. By default, we use an OpenAI 5.2 model, while users can switch to other supported backends, including OpenAI, DeepSeek, Anthropic, Google, and Ollama. Based on extensive internal testing, the OpenAI 5.2 model provides the best overall performance and thus serves as the default choice, achieving a success rate above $90\%$. Here, the success rate is defined as the probability that a valid geometry is generated from a given prompt. If a generation attempt fails, users can simply start a new conversation (by clicking \texttt{new chat}) and retry.
	
	It is worth noting that the prompt should specify not only the geometry but also the locations of boundary conditions. If boundary-condition locations are not explicitly provided, LM-DEM will automatically select them. In general, a single LLM call is often sufficient for 2D cases, whereas 3D geometry generation typically requires multiple refinement iterations. In LM-DEM, such an iteration refers to feeding error messages from a failed attempt back to the LLM, allowing it to revise the previous output.
	
	In addition to LLM-based generation, users may import their own Gmsh mesh files (typically in \texttt{.msh} format). After an LLM-generated geometry is obtained, users can further refine it in two ways: (i) directly editing the generated \texttt{.geo} script, or (ii) opening the geometry in the Gmsh GUI for interactive modifications. Currently, the web version \url{https://ai4m.llmdem.com} does not support Gmsh GUI-based editing; this feature is available only in the locally deployed version.
	
\begin{figure}
		\centering
		\includegraphics[scale=0.4]{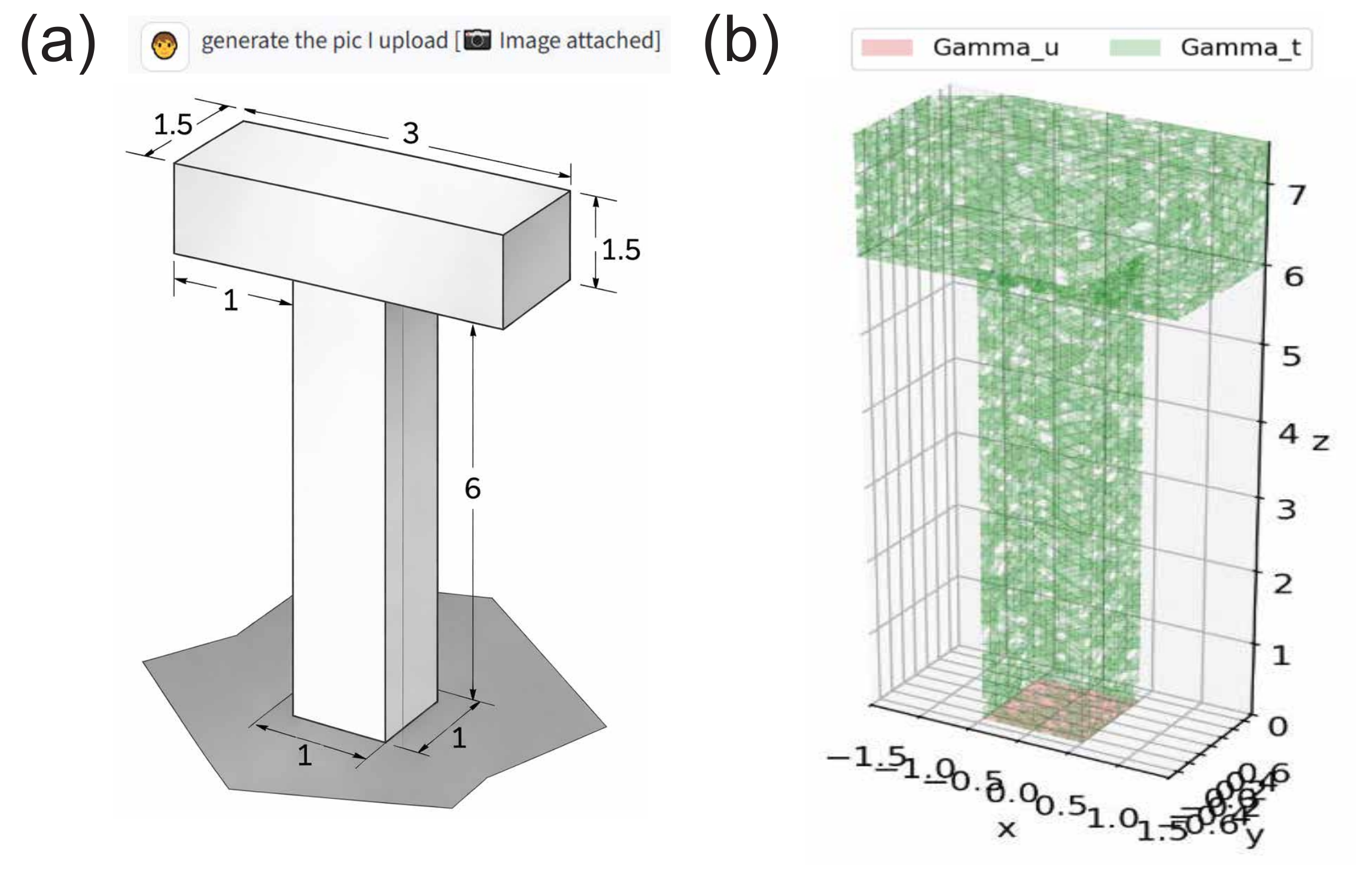}
		\caption{LLM-assisted geometry modeling in LM-DEM with image input: (a) the input image (e.g., a T-shaped structure); (b) the generated geometry and mesh. LM-DEM supports multimodal input: in addition to text prompts, users can upload an image and the LLM infers the geometry and boundary-condition layout, producing a \texttt{.geo} file for Gmsh.}
		\label{fig:LLM_generate_pic}
	\end{figure}

	LM-DEM also supports multimodal inputs, where users can provide images in addition to text prompts. For instance, by inputting an image of a T-shaped structure, the LLM can generate a corresponding geometry and mesh, as shown in Figure~\ref{fig:LLM_generate_pic}.
	
	\subsection{Solver Configuration}
	
	After the geometry is created, the user specifies the solver parameters. In LM-DEM, the configuration is organized into five components.
	
	\subsubsection{Geometry}
	
	Once the geometry is generated, users may adjust the mesh resolution by specifying a global characteristic length $l_{c}$. The problem dimension (2D or 3D) must be selected before the LLM-based geometry generation, as it directly affects the generated geometry and meshing strategy. Besides LLM-based modeling, users can import manually created \texttt{.msh} files.
	
	LM-DEM enforces consistent naming conventions for physical groups: the domain is denoted by \texttt{Omega}, while the boundary is divided into the Dirichlet boundary \texttt{Gamma\_u} and the Neumann boundary \texttt{Gamma\_t}. Users can further select quadrature rules (e.g., for triangles/quadrilaterals and boundary integration). For derivative evaluation, LM-DEM provides two options: finite-element shape functions and automatic differentiation. We recommend finite-element shape functions, which often deliver higher accuracy and efficiency in practice.
	The comparison between finite-element shape functions (explicit differentiation) and automatic differentiation is detailed in \ref{sec:computational-graph}.
	
	In most cases, users only need to specify the dimension, while the remaining geometry settings can be left at their default values.
	
	\subsubsection{Boundary Conditions}
	
	After identifying the locations of boundary conditions, users specify their values, including Dirichlet constraints and Neumann loads. Analytical expressions can also be provided for spatially varying boundary conditions. LM-DEM supports two approaches to enforce essential (Dirichlet) boundary conditions: the distance method and the penalty method. For the penalty method, users can choose the penalty factor. The distance function corresponds to Eq.~\eqref{eq:distance_function}. However, for certain problems, a naive distance-function enforcement may induce noticeable oscillations near the boundary. To mitigate this issue, we introduce a smooth distance strategy, detailed in \ref{sec:Smooth_distance}.
	
	\subsubsection{Material}
	
	Users specify the governing physics (i.e., the ``material model'') through built-in problem types. In 2D, LM-DEM currently supports Poisson, screened Poisson, linear elasticity, and hyperelasticity. In 3D, LM-DEM supports Poisson and linear elasticity.
	
	Beyond these common computational mechanics problems, users may define a custom energy functional via a dedicated interface by selecting \texttt{custom (user-defined)}, analogous to the \texttt{UMAT} interface in Abaqus. In the custom functional setting, users specify the internal energy and the external work. In most cases, the external work follows the standard inner-product form on the Neumann boundary, e.g., $\int_{\Gamma}\boldsymbol{t}\cdot\boldsymbol{u}\,\mathrm{d}\Gamma$ for linear elasticity, where $\boldsymbol{t}$ is the prescribed traction, $\boldsymbol{u}$ is the displacement field of interest, and $\Gamma$ denotes the corresponding boundary. Therefore, users typically only need to define the internal energy term.
	
	\subsubsection{Network}
	
	LM-DEM provides multiple neural network architectures, including multilayer perceptrons (MLPs) and Kolmogorov--Arnold networks (KAN) \citep{wang2025physics,liu2024kan}. 
	In general, KAN tend to be less efficient than MLPs but can achieve higher accuracy; they are particularly suitable for problems with high-frequency solution features \citep{wang2025physics}.
	Users can specify the number of layers, the width of each layer, and the activation function (e.g., \texttt{tanh}, \texttt{silu}, \texttt{gelu}).
	
	\subsubsection{Training}
	
	LM-DEM allows users to configure training details such as learning rate and the maximum number of epochs. Since DEM may converge early, and overly long training can lead to non-physical fracture-like artifacts \citep{wang2026failure,wang2025towards}, LM-DEM incorporates an early-stopping strategy: once the loss exhibits abnormally large fluctuations after apparent convergence, the optimization is terminated to prevent non-physical behavior.
	
	In our implementation, the DEM neural networks are optimized using SOAP (Shampoo with Adam in the Preconditioner’s eigenbasis) \cite{vyas2024soap}. The code of SOAP is \url{https://github.com/nikhilvyas/SOAP}. SOAP integrates Adam-style moment estimation with a Shampoo-inspired second-order preconditioning scheme, where the preconditioner is updated periodically via eigendecomposition.
	Empirically, SOAP often yields strong performance for PINNs training \cite{wang2025gradient}.
	
	Users may also set the number of training points in the domain and on the boundary, and select whether to use GPU acceleration. At present, the web service \url{https://ai4m.llmdem.com} supports CPU only.
	
	\subsection{Post-processing}
	
	LM-DEM supports both in-app visualization in Streamlit and exporting results to ParaView-compatible files, providing an enhanced post-processing experience.
	
\begin{figure}
		\centering
		\includegraphics[scale=0.4]{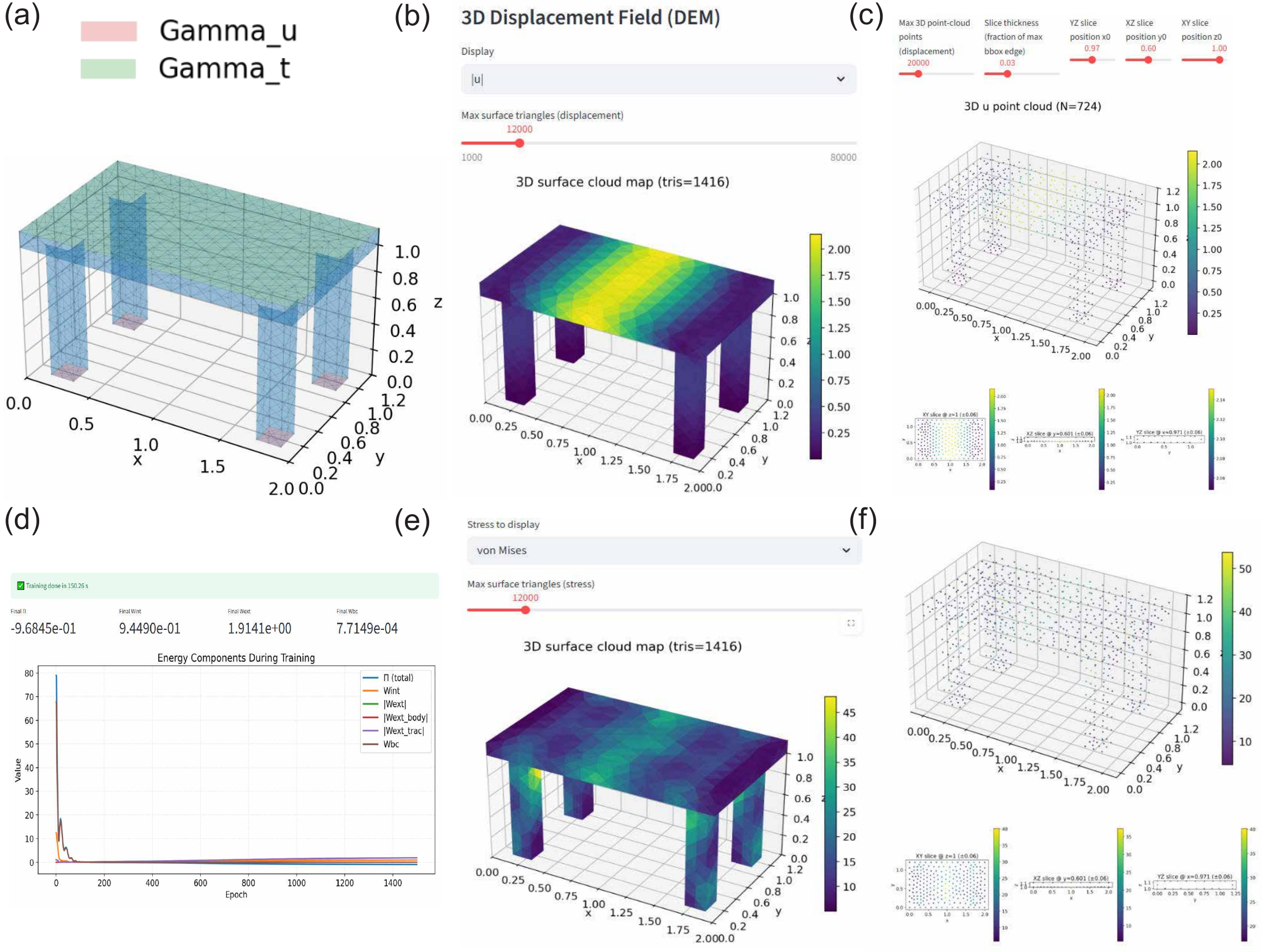}
		\caption{Post-processing in LM-DEM: (a) an LLM-generated 3D geometry with the prompt ``Generate a table; the seating area is the force boundary, and the bottom edges of the four legs are the displacement boundaries.'' (b) 3D displacement contour; (c) displacement point cloud with slicing visualization; (d) evolution of the loss/energy; (e) 3D stress contour; (f) stress point cloud with slicing visualization. LM-DEM provides in-app visualization in the Streamlit interface and supports export to ParaView-compatible \texttt{.vtk} files for further analysis.}
		\label{fig:LM-DEM_post}
	\end{figure}

	In the Streamlit interface, LM-DEM visualizes the evolution of the energy (loss) over epochs as well as the predicted field contours. For 3D problems, LM-DEM additionally provides slicing views and the distribution of training points, as illustrated in Figure~\ref{fig:LM-DEM_post}.
	
	\section{Results}
	\label{sec:Result}
	
	In this section, we present representative results obtained using LM-DEM, including large language model (LLM)-based geometry generation, geometry refinement, numerical solutions computed by the deep energy method and the finite element method, as well as examples involving user-defined energy functionals.
	
	\subsection{Large Language Model-assisted Geometry Modeling}
	
	Geometry generation using large language models is highly intuitive in LM-DEM: users only need to input a textual description of the desired geometry into the dialog box. In the following, we demonstrate both 2D and 3D examples.
	
\begin{figure}
		\centering
		\includegraphics[scale=0.43]{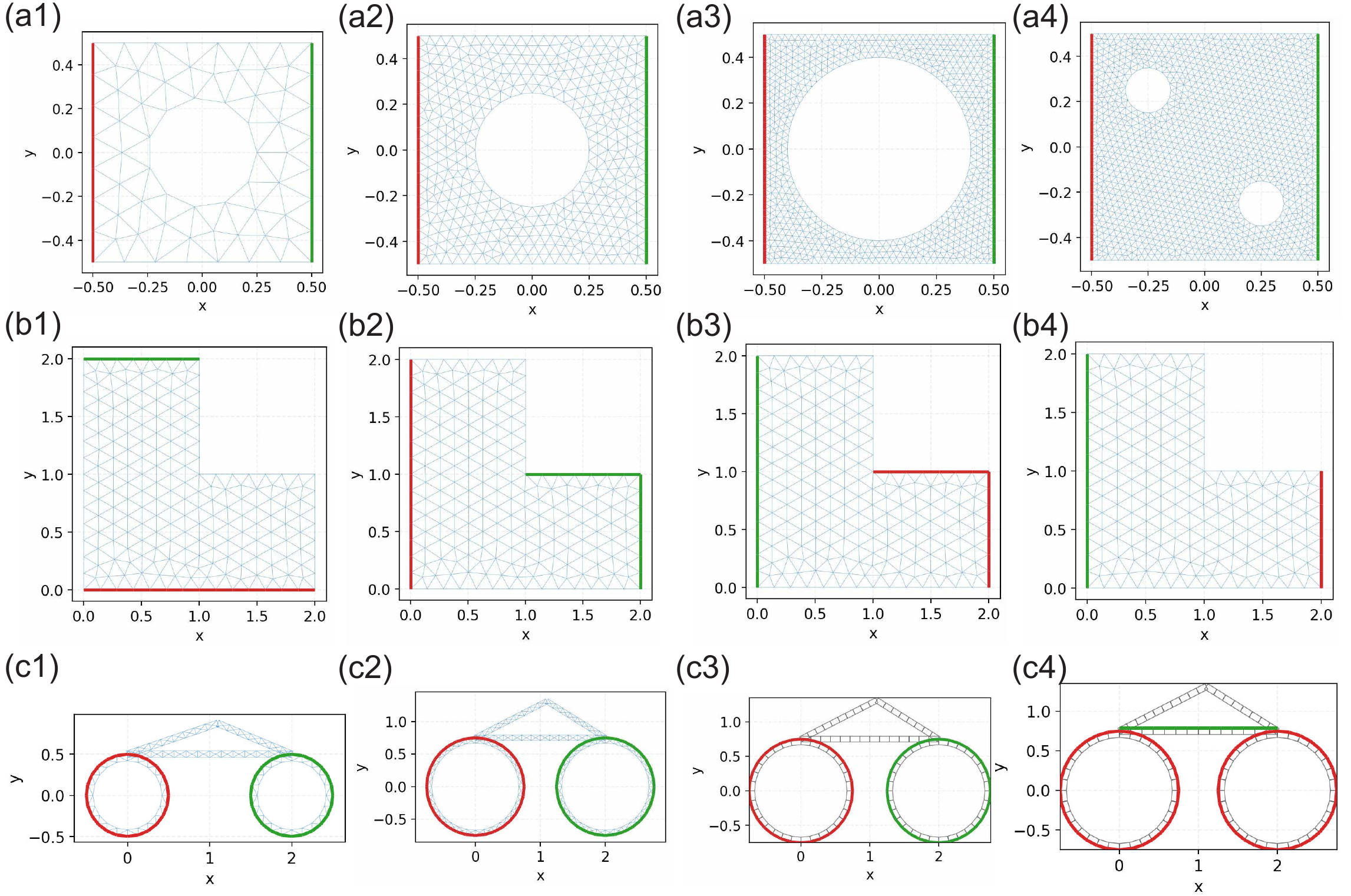}
		\caption{2D geometry examples generated by LM-DEM using LLMs. From top to bottom: a plate with a circular hole, an L-shaped beam, and a bicycle geometry. Each row corresponds to an independent dialog session. From left to right, the columns show successive refinement results based on the prompts listed in \Cref{2D_prompt}. Red regions denote displacement (Dirichlet) boundaries, while green regions indicate force (Neumann) boundaries.}
		\label{fig:LLM_2D}
	\end{figure}

\begin{figure}
		\centering
		\includegraphics[scale=0.42]{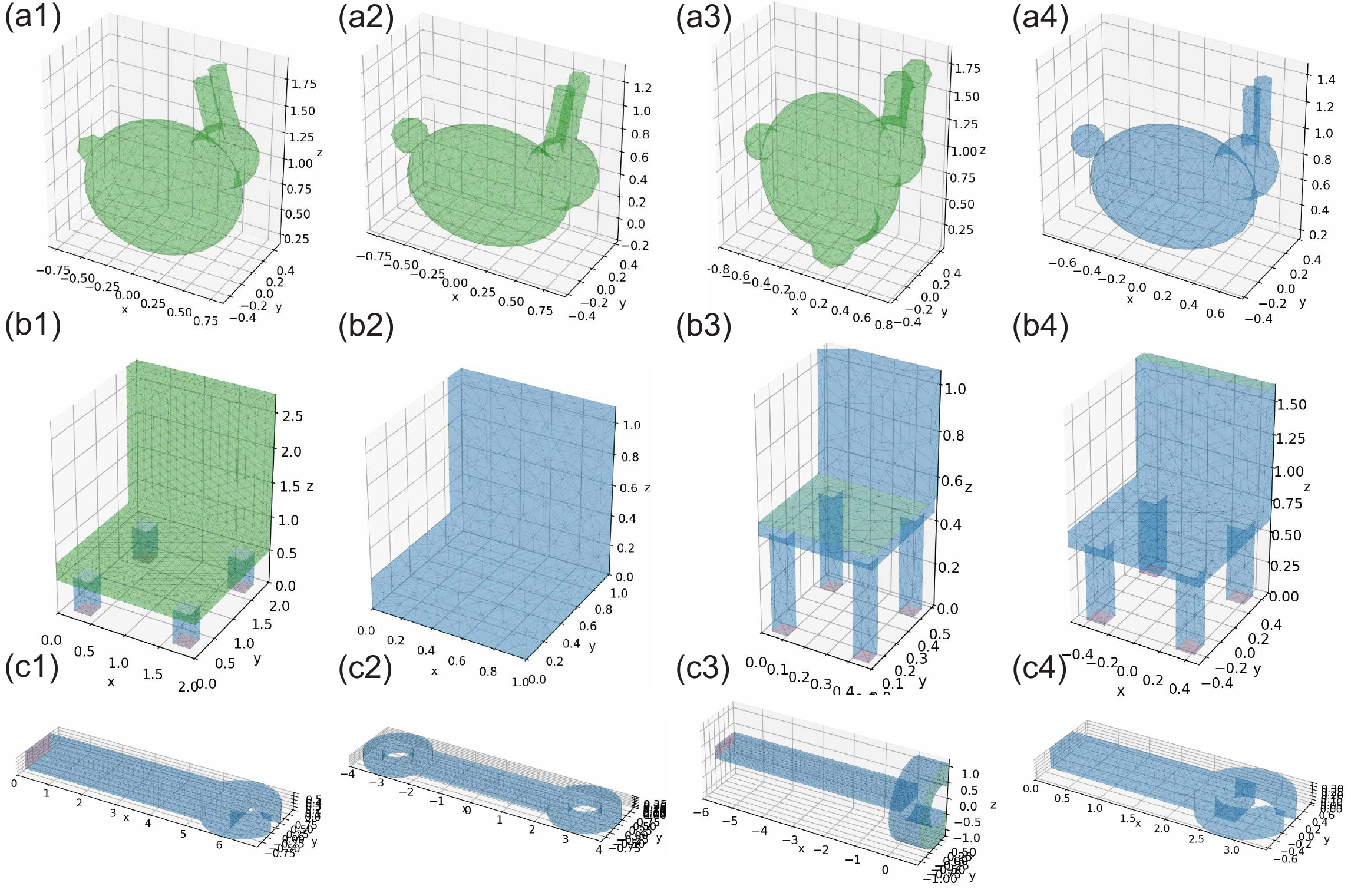}
		\caption{3D geometry examples generated by LM-DEM using LLMs. From top to bottom, the prompts are ``rabbit'', ``chair'', and ``spanner''. Different columns correspond to independent geometry realizations generated from repeated calls to the LLM using the same prompt; the LLM exhibits variability in shape and detail, so users can pick a realization and refine the corresponding \texttt{.geo} file or start a new chat to regenerate.}
		\label{fig:LLM_3D}
	\end{figure}

	Figure~\ref{fig:LLM_2D} shows several 2D geometries generated by the LLM together with their refinement results. It can be observed that the LLM is able to accurately interpret the given prompts and construct the corresponding geometries. We then consider a three-dimensional case. After providing a brief description of the problem, the computational geometry is automatically generated by repeatedly invoking the LLM. The generated geometries are subsequently visualized and analyzed to evaluate the LLM’s capability in modeling complex 3D shapes. Figure~\ref{fig:LLM_3D} presents multiple geometry realizations obtained from repeated calls to the LLM using the same prompt. It is evident that the LLM exhibits strong creativity, producing distinct geometries at each invocation. In principle, even a simple prompt may correspond to infinitely many admissible geometries. Therefore, when geometric dimensions are not strictly specified, users may first generate candidate geometries using the LLM and then manually refine the resulting \texttt{.geo} files. Notably, for every generated geometry, the corresponding \texttt{.geo} file is provided to the user for further modification.
	
	\begin{listing}[htbp]
		\centering
		\fbox{\begin{minipage}{0.92\linewidth}
				\vspace{0.5em}
				\noindent\textbf{Example prompts for 2D geometry generation:}
				\begin{enumerate}
					\item \textbf{Plate with a central circular hole}: 
					\begin{itemize}
						\item Generate a plate with a central circular hole, a $1\times1$ square domain with a hole of radius $0.25$ at the center. The left boundary is subject to displacement constraints, and the right boundary is subjected to force boundary conditions.
						\item Increase elements.
						\item Change the hole radius to $0.4$.
						\item Change to a double-hole configuration: two circular holes with centers at $(-0.25,\,0.25)$ and $(0.25,\,-0.25)$, each with radius $0.1$.
					\end{itemize}
					\item \textbf{L-shaped beam}: 
					\begin{itemize}
						\item Generate a L-shape beam, the bottom boundary is fully fixed, and the top boundary is subjected to force boundary conditions.
						\item The left boundary is subject to displacement (Dirichlet) boundary conditions, and the right boundary is subjected to force (Neumann) boundary conditions.
						\item Swap the boundary conditions.
						\item The displacement (Dirichlet) boundary is only on the right side.
					\end{itemize}
					\item \textbf{Simplified 2D bicycle geometry}: 
					\begin{itemize}
						\item Generate a simplified 2D bicycle geometry, two circular wheels and a triangular frame. The rear and front wheels are circles with centers at $(0,0)$ and $(2,0)$, respectively. The frame consists of three straight beams forming a triangle connecting the top of the rear wheel, the top of the front wheel, and an intermediate seat-tube joint.
						\item Change the tire radius to $1$.
						\item Use quadrilateral elements.
						\item Set the saddle (seat) as a force (Neumann) boundary, and set all tires as displacement (Dirichlet) boundaries.
					\end{itemize}
				\end{enumerate}
				\vspace{0.5em}
		\end{minipage}}
		\caption{Example prompts for 2D geometry generation.}
		\label{2D_prompt}
	\end{listing}
	
	\subsection{Deep Energy Method and Finite Element Solutions}
	
	LM-DEM provides two types of solvers: the deep energy method (DEM) \citep{loss_is_minimum_potential_energy} and the finite element method (FEM). Users may choose either solver depending on their needs. Currently, in 2D, LM-DEM supports Poisson, screened Poisson, linear elasticity, and hyperelasticity problems, while in 3D it supports Poisson and linear elasticity problems. All of these problems can be solved using DEM. However, at present, only linear finite elements are implemented in FEM, which means that 2D hyperelasticity problems cannot yet be solved by FEM within LM-DEM. In addition, for 2D problems, FEM supports both triangular and quadrilateral meshes, whereas for 3D problems only tetrahedral meshes are currently available.
	
\begin{figure}
		\centering
		\includegraphics[scale=0.42]{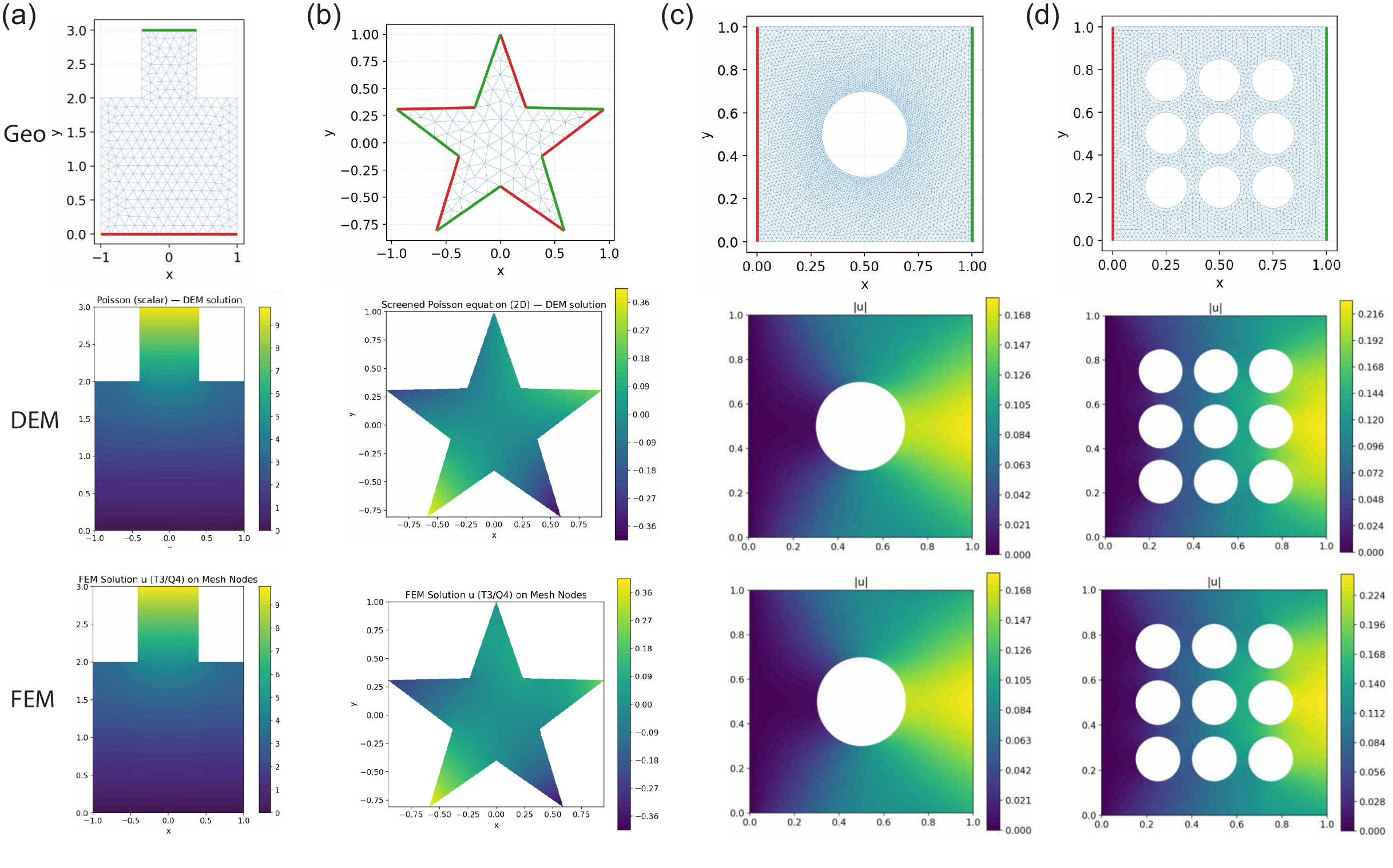}
		\caption{Numerical results obtained by LM-DEM using DEM and FEM for 2D problems. From the first row to the third row: geometry generated from prompts, DEM-predicted displacement magnitude, and FEM-predicted displacement magnitude. Red regions indicate displacement (Dirichlet) boundaries, and green regions indicate force (Neumann) boundaries. From left to right: Poisson problem, screened Poisson equation, and two linear elasticity problems.}
		\label{fig:LM-DEM_2D_DEM_FEM}
	\end{figure}

	We first present 2D results for the Poisson equation, the screened Poisson equation, and linear elasticity. The governing PDEs for these problems are described in \ref{sec:Details-of-PDEs}. For the Poisson problem, the displacement boundary condition is fully fixed, and the traction boundary condition is $t=5$. For the screened Poisson equation, the displacement boundary condition is prescribed as $u=\sin(x)\sin(y)$, and the traction boundary condition is zero. For the two linear elasticity examples, the displacement boundary is fully fixed and the applied traction is $t_x=10$. Figure~\ref{fig:LM-DEM_2D_DEM_FEM} shows the corresponding numerical results obtained by LM-DEM. It is worth emphasizing that the FEM results are not exact solutions either; similar to DEM, FEM is also a numerical solver. Users may flexibly adjust various training and solver settings according to their needs, including network architectures, training strategies, and the enforcement methods for displacement boundary conditions.
	
\begin{figure}
		\centering
		\includegraphics[scale=0.42]{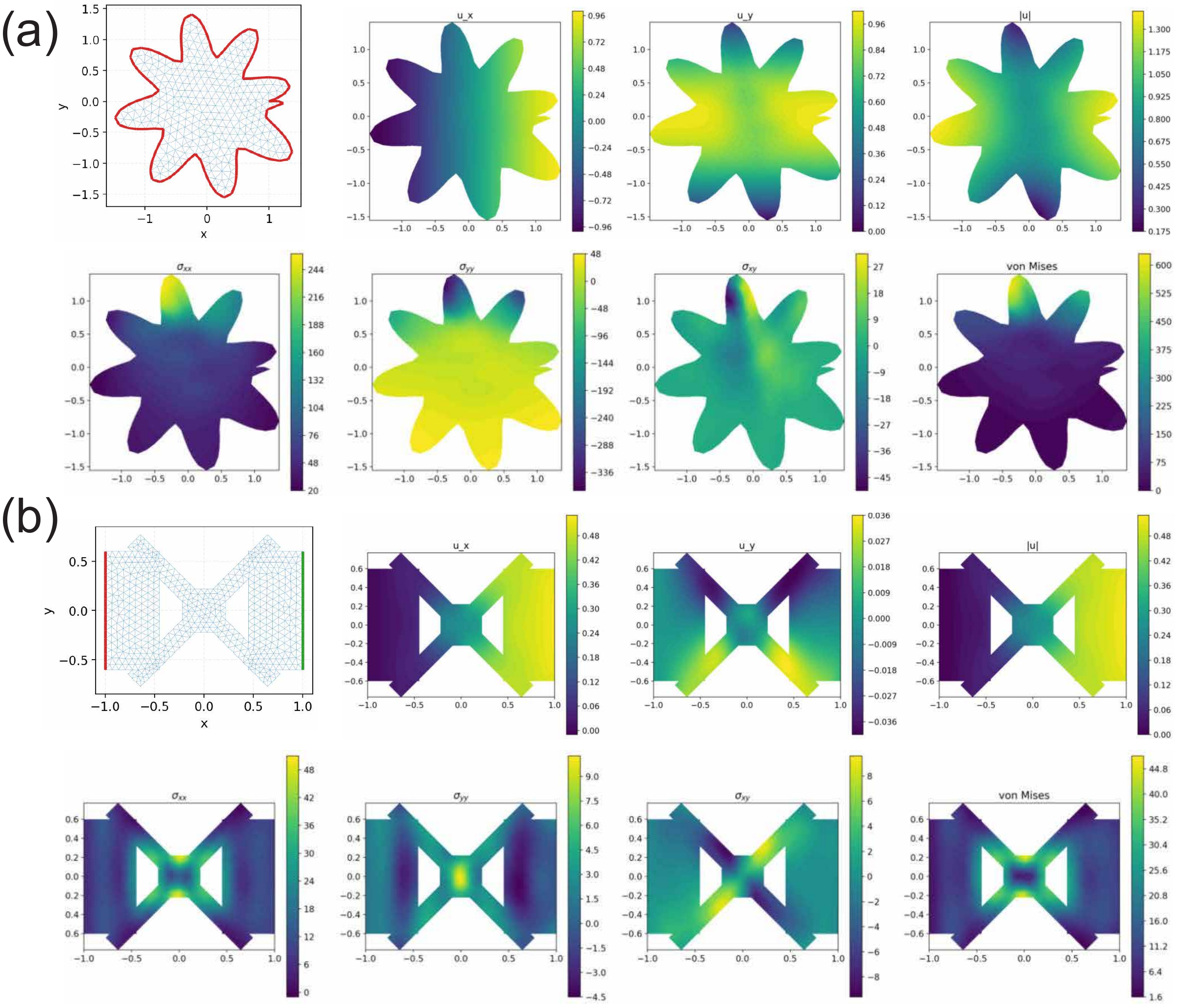}
		\caption{DEM results obtained by LM-DEM for nonlinear problems: (a) spiral petal domain with boundary conditions $u_x=\sin(x)$ and $u_y=\cos(y)$; (b) butterfly load-path structure with the left boundary fully fixed and traction $t_x=10$ on the right boundary. Both are hyperelastic; displacement magnitude and stress contours are shown. Red regions indicate displacement (Dirichlet) boundaries, and green regions indicate force (Neumann) boundaries.}
		\label{fig:LM-DEM_2D_DEM}
	\end{figure}

	Next, we demonstrate the performance of LM-DEM on nonlinear hyperelastic problems. For hyperelasticity, LM-DEM minimizes the total potential energy in a single optimization procedure to obtain the final equilibrium configuration, which is a practical advantage over standard nonlinear FEM workflows that often require load-step control for robust convergence.
	Therefore, for nonlinear problems in the current version of LM-DEM, we provide DEM-based solutions, while FEM is not yet available.
	On the other hand, if one is interested in the intermediate response before the final equilibrium is reached (e.g., to extract a load--deflection curve), the external loading still needs to be applied incrementally and the energy minimization should be performed across successive load steps.
	Figure~\ref{fig:LM-DEM_2D_DEM} presents DEM-predicted displacement and stress contours for a spiral petal domain and a butterfly load-path structure. For the spiral petal domain, the boundary conditions are prescribed as $u_x=\sin(x)$ and $u_y=\cos(y)$. For the butterfly load-path structure, the left boundary is fully fixed, and a traction $t_x=10$ is applied on the right boundary. 
	The prompts used for Figure~\ref{fig:LM-DEM_2D_DEM} are ``spiral petal domain'' and ``butterfly load-path structure''. Note that these prompts are intentionally simple; therefore, repeated queries with the same prompt may generate similar but not identical geometries.

\begin{figure}
		\centering
		\includegraphics[scale=0.40]{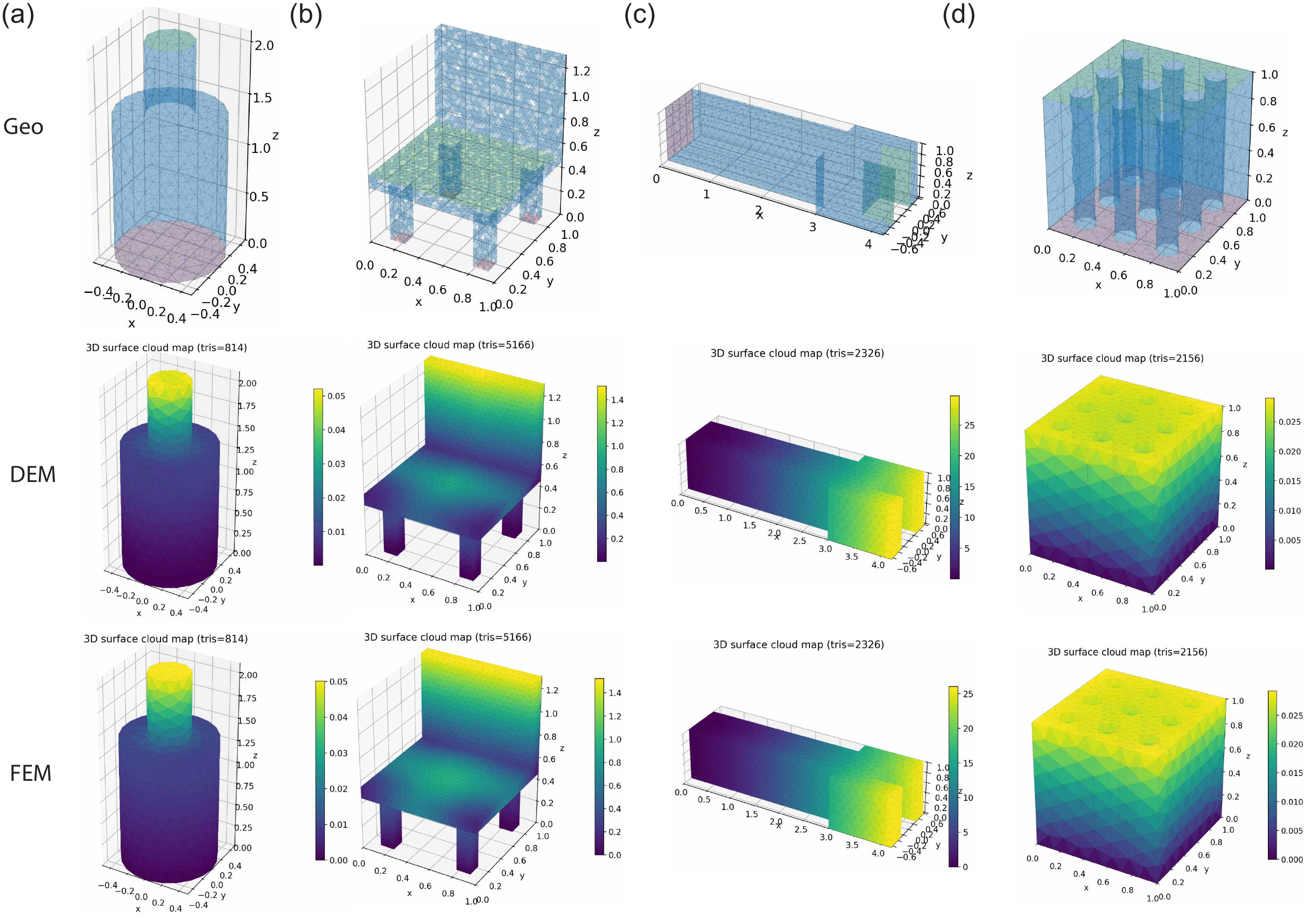}
		\caption{DEM and FEM results obtained by LM-DEM for 3D linear elasticity problems. From the first row to the third row: geometry generated from prompts, DEM-predicted displacement magnitude, and FEM-predicted displacement magnitude. In the first column (bottle), the bottom surface at $z=0$ is fully fixed and a traction $t_z=5$ is applied on the top surface at $z=2$. In the second column (chair), the bottom surface at $z=0$ is fixed and a traction $t_z=-5$ is applied on the seating area at $z=0.5$. In the third column (spanner), the left end at $x=0$ is fixed and tractions $t_y=-5$ are applied on the two clamping regions. In the fourth column (cube with multiple cylindrical holes), the bottom surface at $z=0$ is fixed and a traction $t_z=-3$ is applied on the top surface at $z=1$.}
		\label{fig:LM-DEM_3D_DEM_FEM}
	\end{figure}

	LM-DEM is also capable of solving three-dimensional problems. 
	Figure~\ref{fig:LM-DEM_3D_DEM_FEM} shows DEM and FEM results for various 3D geometries.
	The corresponding prompts are ``bottle'', ``chair'', ``spanner'', and ``cube with multiple cylindrical holes''.
	
	Overall, these examples demonstrate that LM-DEM enables users to conveniently obtain DEM and FEM solutions once the geometry and solver settings are specified. Depending on the application requirements, users can flexibly choose the appropriate solver.
	
	\subsection{User-defined Energy Functionals}
	
	In this subsection, we demonstrate the capability of LM-DEM to handle user-defined energy functionals, in a manner analogous to the \texttt{UMAT} interface in Abaqus. As an illustrative example, we consider hyperelastic constitutive models, including the Neo-Hookean model, the Isihara model \citep{isihara1951statistical}, and the Gent--Thomas model \citep{gent1958forms}. The corresponding constitutive relations are summarized in \ref{sec:hyper_consti}.
	
\begin{figure}
		\centering
		\includegraphics[scale=0.42]{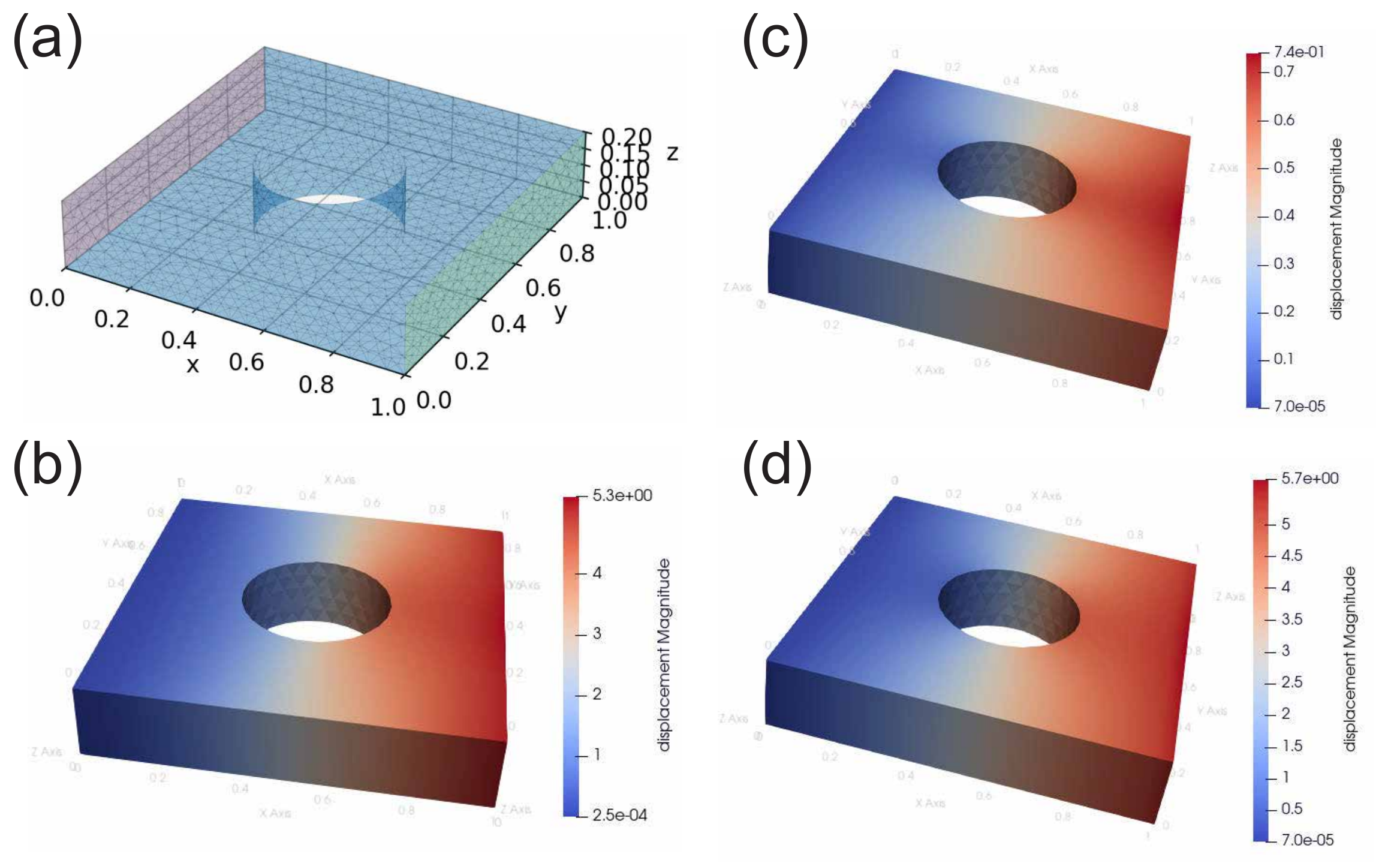}
		\caption{Results obtained by LM-DEM with user-defined constitutive models in 3D: (a) geometry and boundary conditions of a plate with a central hole (fixed at $x=0$, traction $t_x=3$ at $x=1$), where red regions indicate displacement (Dirichlet) boundaries and green regions indicate force (Neumann) boundaries; (b)--(d) displacement magnitude contours for the Neo-Hookean, Isihara, and Gent--Thomas models, respectively. Users define the strain energy density via the custom (user-defined) interface, analogous to the UMAT interface in Abaqus.}
		\label{fig:LM-DEM_user_define}
	\end{figure}

	We generate a three-dimensional plate-with-a-hole geometry, and Figure~\ref{fig:LM-DEM_user_define} shows the displacement magnitude contours obtained using these user-defined constitutive models. The displacement boundary condition is prescribed as fully fixed at $x=0$, while a traction $t_x=3$ is applied at $x=1$.
	
	Note that the copy and paste should be on one line, not multiple lines, in the user-defined energy box.
	
	\section{Conclusion}
	\label{sec:Conclusion}
	
	In this work, we developed an open-source, Streamlit-based software framework for energy-form PINNs, termed \textbf{LM-DEM}. The source code is publicly available at \url{https://github.com/yizheng-wang/LMDEM}, and a web-based version can be accessed at \url{https://ai4m.llmdem.com}. A distinctive feature of LM-DEM is its integration of large language models (LLMs) for geometry modeling: users can generate Gmsh-compatible \texttt{.geo} files directly from natural language descriptions or images. Multiple LLM backends are supported, including OpenAI and DeepSeek. In addition, LM-DEM provides a comprehensive set of modules required for simulation, covering geometry definition, boundary conditions, material models, neural network architectures, and training configurations. Users can further define custom material models through user-defined energy functionals. Once all parameters are specified, solutions can be obtained using the deep energy method, and corresponding finite element solutions can also be generated for comparison. For visualization purposes, LM-DEM supports exporting results in \texttt{.vtk} format to accommodate different post-processing needs.
	
	At present, LM-DEM is limited to quasi-static problems. Future extensions will focus on dynamic problems, such as path-dependent fracture propagation and fluid dynamics. Moreover, the current implementation supports only Dirichlet and Neumann boundary conditions; mixed boundary conditions—commonly encountered in problems such as plates with holes where displacement constraints are applied in one direction and traction conditions in another—are not yet supported and will be addressed in future work. Beyond PINNs, LM-DEM can also be extended to other physics-driven operator-learning frameworks, such as variationally informed neural operators (VINO) \citep{eshaghi2025variational} and the pretrain finite element method (PFEM) \citep{wang2026pretrain}. 
	It is also promising to incorporate agentic workflows as a future direction: by leveraging LLM-based agents, LM-DEM could automate multi-step tasks such as geometry refinement, boundary-condition specification, hyperparameter tuning, and iterative error-driven debugging of the simulation pipeline. 
	Another important direction is to extend the framework to incremental energy functionals that carry history variables, enabling transient or rate-dependent problems (e.g., plasticity or phase-field fracture) within an energy-based formulation. Such an extension would naturally open the door to transfer learning across load/time steps and to adaptive load-step control strategies, thereby improving robustness and reducing training cost.
	Further improvements may include adaptive sampling strategies (e.g., residual- or energy-driven point refinement), fine-tuning or warm-start training for related problem instances, and domain decomposition with appropriate coupling (e.g., Schwarz-type domain decomposition methods) to enhance computational scalability. These extensions are also closely related to parametric DEM formulations, where one learns solution manifolds across families of parameters and combines them with decomposition strategies for improved efficiency. In the current version, both training and numerical integration are performed directly in the physical domain by evaluating the energy functional and its derivatives at sampled/quadrature points.
	Overall, LM-DEM has significant potential for further development, and its continued evolution will benefit from the collective efforts of the open-source community.

	\section*{Declaration of competing interest}
	The authors declare that they have no known competing financial interests or personal relationships that could
	have appeared to influence the work reported in this paper.
	
	\section*{Acknowledgement}
	The study was supported by the Key Project of the National Natural Science Foundation of China (12332005) and scholarship from Bauhaus University in Weimar.

	\section*{CRediT authorship contribution statement}
	
	\textbf{Yizheng Wang}: Conceptualization, Methodology, Formal analysis, Investigation, Data curation, Validation, Visualization, Writing – original draft, Writing – review \& editing. 
	\textbf{Cosmin Anitescu}: Deploy, review.  
	\textbf{Mohammad Sadegh Eshaghi}: Investigation.  
	\textbf{Xiaoying Zhuang}: Supervision, Writing – review \& editing.  
	\textbf{Timon Rabczuk}: Supervision, Writing – review \& editing, Funding acquisition.  
	\textbf{Yinghua Liu}: Supervision, Funding acquisition.

	\appendix

	\section{Comparison of computational graphs for automatic and explicit differentiation}
	\label{sec:computational-graph}
	
	Automatic differentiation (AD) often incurs rapidly increasing computational cost as the derivative order increases, whereas the computational cost of explicit differentiation remains essentially independent of the derivative order.
	
	Therefore, LM-DEM adopts explicit differentiation (based on finite element shape functions) as the default strategy for computing spatial derivatives, rather than relying on AD. In what follows, we explain why explicit differentiation can be significantly more efficient than AD.
	
	To illustrate this difference, we consider the following PDE:
	\begin{equation}
		\nabla \cdot \big[\nabla T(x,y) \big] = f(x,y),
		\label{eq:computational_graph}
	\end{equation}
	where spatial derivatives of $T(\boldsymbol{x})$ are required in order to evaluate the loss of DEM.
	
	\begin{figure}[t]
		\centering
		\includegraphics[scale=0.52]{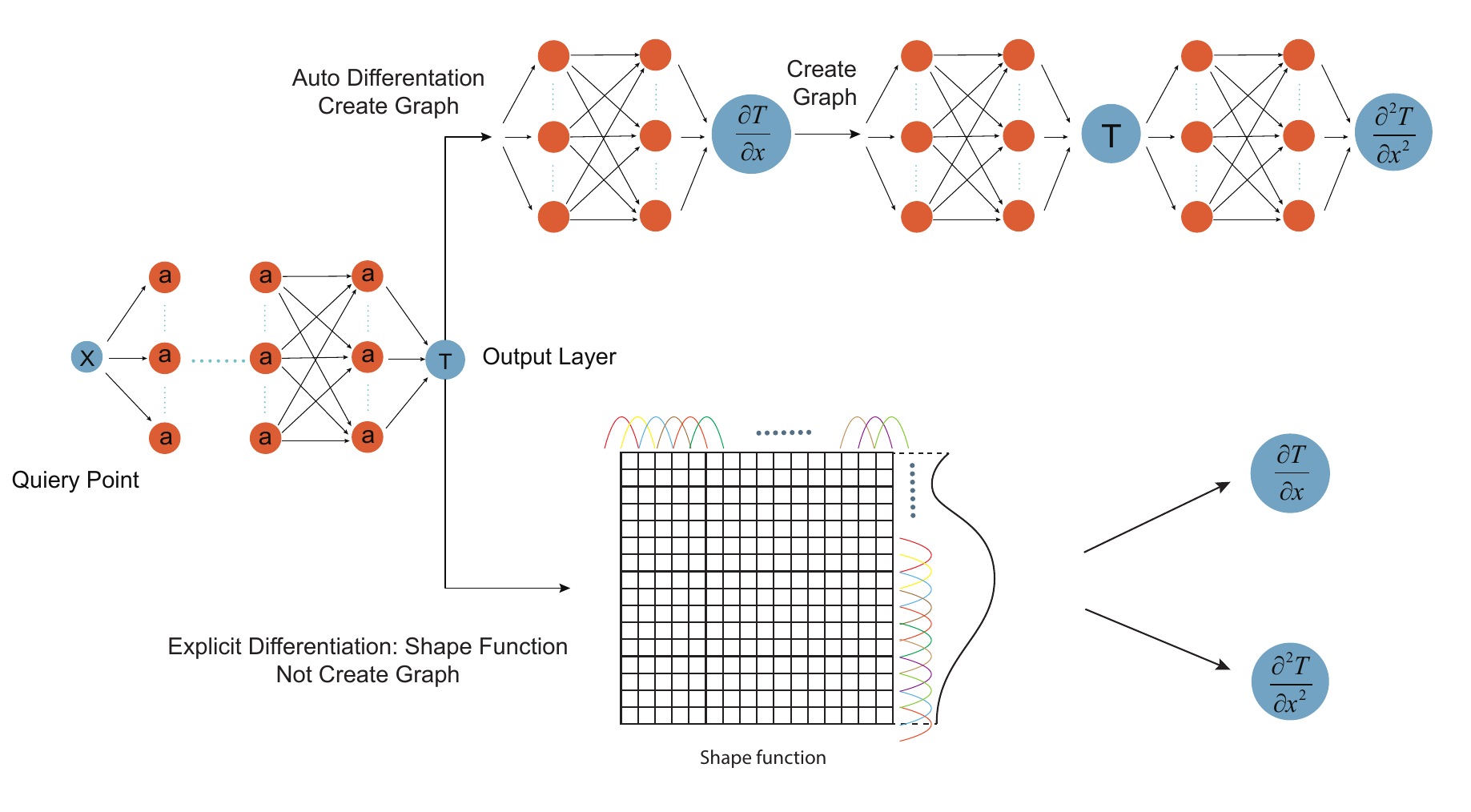}
		\caption{Comparison of computational graphs for automatic differentiation and explicit differentiation.}
		\label{fig:Computation-graph}
	\end{figure}
	
	Below, we analyze the computational complexity of automatic differentiation and explicit differentiation, respectively.
	
	The upper part of \Cref{fig:Computation-graph} illustrates the computational graph associated with automatic differentiation. When constructing the derivative via AD, new computational graphs are generated recursively through repeated application of the chain rule. Let $C_T$ denote the number of basic computational steps required to evaluate $\boldsymbol{x} \mapsto T(\boldsymbol{x})$. Then:
	\begin{itemize}
		\item The computational graph for $\boldsymbol{x} \mapsto \partial T / \partial \boldsymbol{x}$ requires approximately $2C_T$ steps;
		\item The computational graph for $\boldsymbol{x} \mapsto \partial^2 T / \partial \boldsymbol{x}^2$ requires approximately $4C_T$ steps.
	\end{itemize}
	More generally, if the highest derivative order appearing in the PDE is $M$, the total computational complexity scales as
	$
	\mathcal{O}(2^M C_T),
	$
	This exponential growth with respect to the derivative order makes automatic differentiation increasingly inefficient for PDE-based training, especially when higher-order derivatives are involved.
	
	In contrast, when using explicit differentiation by finite element shape functions as shown in the lower part of \Cref{fig:Computation-graph}, spatial derivatives can be derived analytically in closed form, regardless of the derivative order. Once the shape-function representation is established, derivatives of arbitrary order can be evaluated explicitly without introducing additional computational graphs.
	
	As a result, the computational cost of explicit differentiation does not grow with $M$. 
	Although explicit differentiation is typically more efficient (and often more accurate) than AD, it usually requires a mesh in order to define the shape functions. In this sense, DEM with shape-function-based explicit differentiation is no longer strictly meshfree. By contrast, when AD is employed, derivatives can be evaluated at arbitrary points and the energy integrals can be approximated via standard Monte Carlo sampling, yielding a meshfree DEM formulation.

	\section{Smooth Distance Function}
	\label{sec:Smooth_distance}
	
\begin{figure}
		\centering
		\includegraphics[scale=0.5]{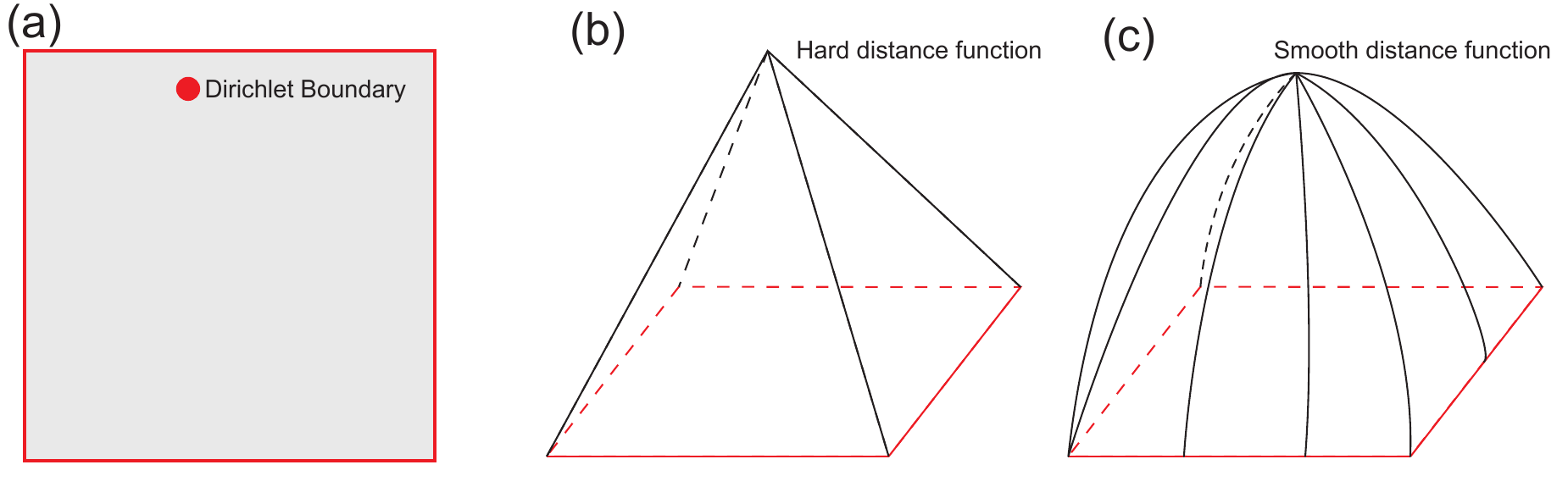}
		\caption{Distance functions: (a) a square domain with all boundaries prescribed as Dirichlet boundaries; (b) schematic illustration of the conventional hard distance function, whose derivative is non-smooth along the black line (where the closest boundary segment switches), leading to gradient discontinuities that can harm energy evaluation in DEM; (c) schematic illustration of the smooth distance function, which is differentiable everywhere and avoids such discontinuities, improving numerical stability; LM-DEM uses the smooth distance combined with a penalty term by default.}
		\label{fig:distance_function}
	\end{figure}

	In the deep energy method (DEM), essential (Dirichlet) boundary conditions must typically be satisfied a priori, which is commonly achieved through the use of a distance function, as defined in Eq.~\eqref{eq:distance_function}. However, for geometries with complex boundaries, the conventional (hard) distance function may introduce non-smoothness at locations where the closest boundary segment switches, as illustrated in Figure~\ref{fig:distance_function}(b). This non-smoothness leads to discontinuities in the gradient $\nabla D$, which in turn may cause inaccuracies in energy evaluation. Since the energy functional usually involves derivatives of the field of interest, differentiation of the distance function can significantly affect the numerical stability of DEM.
	
	To improve the robustness of energy computation, we introduce a smooth distance function, as shown in Figure~\ref{fig:distance_function}(c).
	
	The smooth distance function is defined as
	\begin{equation}
		D_{s}(\boldsymbol{x})
		=
		-\tau \log\!\left[
		\frac{1}{M}
		\sum_{i=1}^{M}
		\exp\!\left(-\frac{d_{i}^{2}}{\tau}\right)
		\right],
	\end{equation}
	where $d_{i}$ denotes the shortest distance from $\boldsymbol{x}$ to the $i$-th segment of the Dirichlet boundary. For complex boundaries, the Dirichlet boundary is partitioned into $M$ segments. The parameter $\tau$ is a smoothing coefficient, which is set to $0.001$ by default in LM-DEM. It can be shown that as $\tau \to 0^{+}$, the smooth distance function $D_{s}(\boldsymbol{x})$ converges to the hard distance function.
	
	The smooth distance function is globally differentiable, as it avoids the non-smooth $\min$ operator inherent in the hard distance function. Nevertheless, a side effect of the smooth distance function is that its value is not exactly zero on the Dirichlet boundary. To address this issue, LM-DEM incorporates an additional penalty term. Consequently, the default strategy in LM-DEM combines the smooth distance function with a penalty method to enforce Dirichlet boundary conditions in a stable and accurate manner.
	
	\section{Details of PDEs}
	\label{sec:Details-of-PDEs}
	
	In this section, we summarize the partial differential equations (PDEs) that are currently built into LM-DEM.
	
	\subsection{Poisson Equation}
	
	The strong form of the Poisson equation is given by
	\begin{equation}
		\begin{cases}
			-\Delta u = f, & \boldsymbol{x}\in\Omega,\\
			u = \bar{u}, & \boldsymbol{x}\in\Gamma_{u},\\
			\dfrac{\partial u}{\partial n} = g, & \boldsymbol{x}\in\Gamma_{t},
		\end{cases}
	\end{equation}
	where $\Gamma_{u}$ and $\Gamma_{t}$ denote the Dirichlet and Neumann boundaries, respectively.
	
	The corresponding energy functional reads
	\begin{align}
		\Pi(u)
		&=
		\int_{\Omega}\frac{1}{2}\nabla u \cdot \nabla u\,\mathrm{d}\Omega
		-
		\int_{\Omega} f u\,\mathrm{d}\Omega
		-
		\int_{\Gamma^{t}} g u\,\mathrm{d}\Gamma,
		\\
		\text{s.t.}\quad
		& u=\bar{u}, \quad \boldsymbol{x}\in\Gamma_{u}.
		\nonumber
	\end{align}
	
	\subsection{Screened Poisson Equation}
	
	The strong form of the screened Poisson equation is
	\begin{equation}
		\begin{cases}
			-\Delta u + k^{2} u = f, & \boldsymbol{x}\in\Omega,\\
			u = \bar{u}, & \boldsymbol{x}\in\Gamma_{u},\\
			\dfrac{\partial u}{\partial n} = g, & \boldsymbol{x}\in\Gamma_{t}.
		\end{cases}
	\end{equation}
	
	The corresponding energy functional is given by
	\begin{align}
		\Pi(u)
		&=
		\int_{\Omega}\frac{1}{2}\nabla u \cdot \nabla u\,\mathrm{d}\Omega
		+
		\int_{\Omega}\frac{1}{2}k^{2}u^{2}\,\mathrm{d}\Omega
		-
		\int_{\Omega} f u\,\mathrm{d}\Omega
		-
		\int_{\Gamma^{t}} g u\,\mathrm{d}\Gamma,
		\\
		\text{s.t.}\quad
		& u=\bar{u}, \quad \boldsymbol{x}\in\Gamma_{u}.
		\nonumber
	\end{align}
	
	\subsection{Linear Elasticity}
	
	Linear elasticity is one of the most classical benchmark problems in solid mechanics. Its strong form is
	\begin{equation}
		\begin{aligned}
			\sigma_{ij,j} + f_{i} &= 0,\\
			\sigma_{ij} &=
			\frac{E}{1+\nu}\varepsilon_{ij}
			+
			\frac{E\nu}{(1+\nu)(1-2\nu)}\varepsilon_{kk}\delta_{ij},\\
			\varepsilon_{ij} &= \frac{1}{2}(u_{i,j}+u_{j,i}),\\
			\boldsymbol{u} &= \bar{\boldsymbol{u}}, \quad \boldsymbol{x}\in\Gamma_{\boldsymbol{u}},\\
			\sigma_{ij}n_{j} &= \bar{t}_{i}, \quad \boldsymbol{x}\in\Gamma_{\boldsymbol{t}},
		\end{aligned}
	\end{equation}
	where $\boldsymbol{\sigma}$ and $\boldsymbol{\varepsilon}$ denote the stress and strain tensors, respectively, and $E$ and $\nu$ are the Young’s modulus and Poisson’s ratio.
	
	The corresponding energy formulation is
	\begin{equation}
		\begin{aligned}
			\Pi(\boldsymbol{u})
			&=
			\int_{\Omega}\big(\Psi - \boldsymbol{f}\cdot\boldsymbol{u}\big)\,\mathrm{d}\Omega
			-
			\int_{\Gamma^{\boldsymbol{t}}}\bar{\boldsymbol{t}}\cdot\boldsymbol{u}\,\mathrm{d}\Gamma,\\
			\Psi &= \frac{1}{2}\sigma_{ij}\varepsilon_{ij},\\
			\text{s.t.}\quad
			& \boldsymbol{u}=\bar{\boldsymbol{u}}, \quad \boldsymbol{x}\in\Gamma_{\boldsymbol{u}},
		\end{aligned}
	\end{equation}
	where $\Psi$ denotes the strain energy density.
	
	\subsection{Hyperelasticity}
	
	Hyperelastic problems are nonlinear with respect to both geometry and material behavior. The strong form is
	\begin{equation}
		\begin{aligned}
			P_{iI,I} + f_{i} &= 0,\\
			\boldsymbol{P} &= \frac{\partial \Psi}{\partial \boldsymbol{F}},\\
			\boldsymbol{F} &= \frac{\partial \boldsymbol{x}}{\partial \boldsymbol{X}},\\
			\boldsymbol{u} &= \bar{\boldsymbol{u}}, \quad \boldsymbol{x}\in\Gamma_{\boldsymbol{u}},\\
			P_{iI}n_{I} &= \bar{t}_{i}, \quad \boldsymbol{x}\in\Gamma_{\boldsymbol{t}},
		\end{aligned}
	\end{equation}
	where $\Psi$ is the strain energy density, $\boldsymbol{P}$ is the first Piola--Kirchhoff stress tensor, and $\boldsymbol{F}$ is the deformation gradient. Here, $\boldsymbol{x}$ and $\boldsymbol{X}$ denote the spatial and material coordinates, respectively, related by $\boldsymbol{x}=\boldsymbol{X}+\boldsymbol{u}$.
	
	The corresponding energy functional is
	\begin{equation}
		\begin{aligned}
			\mathcal{L}
			&=
			\int_{\Omega}\big(\Psi - \boldsymbol{f}\cdot\boldsymbol{u}\big)\,\mathrm{d}\Omega
			-
			\int_{\Gamma^{\boldsymbol{t}}}\bar{\boldsymbol{t}}\cdot\boldsymbol{u}\,\mathrm{d}\Gamma,\\
			\text{s.t.}\quad
			& \boldsymbol{u}=\bar{\boldsymbol{u}}, \quad \boldsymbol{x}\in\Gamma_{\boldsymbol{u}}.
		\end{aligned}
	\end{equation}
	
	For a Neo-Hookean material, the strain energy density is given by
	\begin{equation}
		\begin{aligned}
			\Psi
			&=
			\frac{1}{2}\lambda (\ln J)^{2}
			-
			\mu \ln J
			+
			\frac{1}{2}\mu (I_{1}-2),\\
			J &= \det(\boldsymbol{F}),\\
			I_{1} &= \mathrm{tr}(\boldsymbol{C}),\\
			\boldsymbol{C} &= \boldsymbol{F}^{\mathsf{T}}\boldsymbol{F},
		\end{aligned}
	\end{equation}
	where $\lambda$ and $\mu$ are the Lamé parameters,
	\begin{equation}
		\lambda = \frac{\nu E}{(1+\nu)(1-2\nu)},
		\qquad
		\mu = \frac{E}{2(1+\nu)},
	\end{equation}
	and $\boldsymbol{C}$ denotes the right Cauchy--Green deformation tensor.
	
	\section{Hyperelastic Constitutive Models}
	\label{sec:hyper_consti}
	
	In this section, we introduce the hyperelastic constitutive models used in LM-DEM, including the Neo-Hookean model, the Isihara (IH) model \citep{isihara1951statistical}, and the Gent--Thomas (GT) model \citep{gent1958forms}.
	
	\subsection{Neo-Hookean Model}
	
	The Neo-Hookean model is one of the most fundamental and widely used hyperelastic constitutive models. Its main advantages are the small number of material parameters and strong numerical stability; however, it tends to exhibit excessive softening under large deformations. The strain energy density function of the Neo-Hookean model is given by
	\begin{equation}
		\Psi
		=
		a_{0}(\tilde{I}_{1}-3)
		+
		a_{1}(J-1)^{2},
	\end{equation}
	where $a_{0}$ and $a_{1}$ are material parameters, which are set to $0.5$ and $1.5$, respectively, in this work \citep{thakolkaran2025can}. The quantity $\tilde{I}_{1}$ denotes the first isochoric invariant, defined as
	\begin{equation}
		\tilde{I}_{1}
		=
		\mathrm{tr}(\tilde{\boldsymbol{C}})
		=
		J^{-\frac{2}{3}}\mathrm{tr}(\boldsymbol{C}),
	\end{equation}
	where $J=\det(\boldsymbol{F})$ is the determinant of the deformation gradient $\boldsymbol{F}$, and $\boldsymbol{C}$ is the right Cauchy--Green deformation tensor, defined as $\boldsymbol{C}=\boldsymbol{F}^{\mathsf{T}}\boldsymbol{F}$. The isochoric counterparts are given by $\tilde{\boldsymbol{C}}=\tilde{\boldsymbol{F}}^{\mathsf{T}}\tilde{\boldsymbol{F}}$ and $\tilde{\boldsymbol{F}}=J^{-1/3}\boldsymbol{F}$.
	
	In LM-DEM, the Neo-Hookean strain energy density can be directly specified using symbolic expressions. For clarity, we list below the corresponding expression implemented in LM-DEM:
	\begin{verbatim}
		0.5*(((1+ux)*(1+vy)*(1+wz) + uy*vz*wx + uz*vx*wy
		- uz*(1+vy)*(wx) - uy*vx*(1+wz) - (1+ux)*vz*wy)**(-2/3)
		*((1+ux)**2 + vx**2 + wx**2 + uy**2 + (1+vy)**2 + wy**2
		+ uz**2 + vz**2 + (1+wz)**2) - 3)
		+ 1.5*(((1+ux)*(1+vy)*(1+wz) + uy*vz*wx + uz*vx*wy
		- uz*(1+vy)*(wx) - uy*vx*(1+wz) - (1+ux)*vz*wy)-1)**2
	\end{verbatim}
	Note that the copy and paste should be on one line, not multiple lines, in the user-defined energy box.

	\subsection{Isihara (IH) Model}
	
	The Isihara model is generally stiffer than the Neo-Hookean model and is particularly effective in capturing moderate-to-large shear deformations. Its strain energy density function is expressed as
	\begin{equation}
		\Psi
		=
		a_{0}(\tilde{I}_{1}-3)
		+
		a_{1}(\tilde{I}_{2}-3)
		+
		a_{2}(\tilde{I}_{1}-3)^{2}
		+
		a_{3}(J-1)^{2},
	\end{equation}
	where $a_{0}$, $a_{1}$, $a_{2}$, and $a_{3}$ are material parameters. In this work, they are chosen as $0.5$, $1.0$, $1.0$, and $1.5$, respectively. The second isochoric invariant $\tilde{I}_{2}$ is defined as
	\begin{equation}
		\begin{aligned}
			\tilde{I}_{2}
			&=
			\frac{1}{2}\big[\mathrm{tr}(\tilde{\boldsymbol{C}})^{2}
			-
			\mathrm{tr}(\tilde{\boldsymbol{C}}^{2})\big]
			=
			J^{-\frac{4}{3}} I_{2},\\
			I_{2}
			&=
			\frac{1}{2}\big[\mathrm{tr}(\boldsymbol{C})^{2}
			-
			\mathrm{tr}(\boldsymbol{C}^{2})\big].
		\end{aligned}
	\end{equation}
	
	The corresponding symbolic expression used in LM-DEM for the Isihara model is listed below:
	\begin{verbatim}
		0.5*(((1+ux)*(1+vy)*(1+wz) + uy*vz*wx + uz*vx*wy - uz*(1+vy)*(wx)
		- uy*vx*(1+wz) - (1+ux)*vz*wy)**(-2/3)*((1+ux)**2 + vx**2 + wx**2
		+ uy**2 + (1+vy)**2 + wy**2 + uz**2 + vz**2 + (1+wz)**2) - 3) 
		+ (((1+ux)*(1+vy)*(1+wz) + uy*vz*wx + uz*vx*wy - uz*(1+vy)*(wx) 
		- uy*vx*(1+wz) - (1+ux)*vz*wy)**(-4/3)*0.5*(((1+ux)**2 + vx**2 
		+ wx**2 + uy**2 + (1+vy)**2 + wy**2 + uz**2 + vz**2 + (1+wz)**2)**2
		-(((1+ux)**2+vx**2+wx**2)**2+(uy**2+(1+vy)**2+wy**2)**2+(uz**2+vz**2
		+(1+wz)**2)**2+2*((1+ux)*uy+vx*(1+vy)+wx*wy)**2+2*((1+ux)*uz+vx*vz+wx
		*(1+wz))**2+2*(uy*uz+(1+vy)*vz+wy*(1+wz))**2))-3) + (((1+ux)*(1+vy)
		*(1+wz) + uy*vz*wx + uz*vx*wy - uz*(1+vy)*(wx) - uy*vx*(1+wz) 
		- (1+ux)*vz*wy)**(-2/3)*((1+ux)**2 + vx**2 + wx**2 + uy**2 + (1+vy)**2 
		+ wy**2 + uz**2 + vz**2 + (1+wz)**2)-3)**2 + 1.5*( ((1+ux)*(1+vy)*(1+wz) 
		+ uy*vz*wx + uz*vx*wy - uz*(1+vy)*(wx) - uy*vx*(1+wz) - (1+ux)*vz*wy)-1)**2
	\end{verbatim}
	where the full expression follows directly from the definitions of $\tilde{I}_{1}$, $\tilde{I}_{2}$, and $J$ as implemented in LM-DEM.
	
	Note that the copy and paste should be on one line, not multiple lines, in the user-defined energy box.
	
	\subsection{Gent--Thomas (GT) Model}
	
	The Gent--Thomas model is designed to describe material behavior approaching limiting stretch. Its strain energy density function is given by
	\begin{equation}
		\Psi
		=
		a_{0}(\tilde{I}_{1}-3)
		+
		a_{1}\log\!\left(\frac{\tilde{I}_{2}}{3}\right)
		+
		a_{2}(J-1)^{2},
	\end{equation}
	where $a_{0}$, $a_{1}$, and $a_{2}$ are material parameters. In this work, these parameters are set to $0.5$, $1.0$, and $1.5$, respectively.
	
	Similarly, the Gent--Thomas model can be specified in LM-DEM using a symbolic expression. The symbolic expression used in LM-DEM for the Gent--Thomas model is listed below:
	\begin{verbatim}
		0.5*(((1+ux)*(1+vy)*(1+wz) + uy*vz*wx + uz*vx*wy - uz*(1+vy)*(wx)
		- uy*vx*(1+wz) - (1+ux)*vz*wy)**(-2/3)*((1+ux)**2 + vx**2 + wx**2
		+ uy**2 + (1+vy)**2 + wy**2 + uz**2 + vz**2 + (1+wz)**2) - 3)
		+ log(((1+ux)*(1+vy)*(1+wz) + uy*vz*wx + uz*vx*wy - uz*(1+vy)*(wx)
		- uy*vx*(1+wz) - (1+ux)*vz*wy)**(-4/3)*0.5*(((1+ux)**2 + vx**2
		+ wx**2 + uy**2 + (1+vy)**2 + wy**2 + uz**2 + vz**2 + (1+wz)**2)**2
		-(((1+ux)**2+vx**2+wx**2)**2+(uy**2+(1+vy)**2+wy**2)**2+(uz**2+vz**2
		+(1+wz)**2)**2+2*((1+ux)*uy+vx*(1+vy)+wx*wy)**2+2*((1+ux)*uz+vx*vz
		+wx*(1+wz))**2+2*(uy*uz+(1+vy)*vz+wy*(1+wz))**2))/3)
		+ 1.5*((1+ux)*(1+vy)*(1+wz) + uy*vz*wx + uz*vx*wy - uz*(1+vy)*(wx)
		- uy*vx*(1+wz) - (1+ux)*vz*wy-1)**2
	\end{verbatim}
	The corresponding implementation follows directly from the above strain energy density formulation and the definitions of $\tilde{I}_{1}$, $\tilde{I}_{2}$, and $J$.
	
	Note that the copy and paste should be on one line, not multiple lines, in the user-defined energy box.
	
	\section{Supplementary code}
	The code of LM-DEM is available at \url{https://github.com/yizheng-wang/LMDEM}. The website of LM-DEM is \url{https://ai4m.llmdem.com}. 
	The links to the instructional video are: \url{https://space.bilibili.com/327201044/lists?sid=7309959&spm_id_from=333.788.0.0} (Chinese) and \url{https://www.youtube.com/channel/UCu0MPs07oYViB5gR4sN46Lg} (English).

	\bibliographystyle{elsarticle-num}
	\addcontentsline{toc}{section}{\refname}\bibliography{reference.bib}
	
\end{document}